\documentclass[11pt,twoside,a4paper]{article}

\usepackage{amsmath}
\usepackage{amsfonts}
\usepackage{cite}

\newcommand{\Ent}[1]{[\mkern - 2.5 mu [#1] \mkern - 2.5 mu ]}  

\begin{document}

\thispagestyle{empty}

\title{Prediction Properties of Aitken's Iterated $\Delta^2$ Process,\\ 
of Wynn's Epsilon Algorithm,\\ 
and of Brezinski's Iterated Theta Algorithm}
\author{Ernst Joachim Weniger \\
Institut f\"{u}r Physikalische und Theoretische Chemie \\
Universit\"{a}t Regensburg \\
D-93040 Regensburg \\
Germany \\
joachim.weniger@chemie.uni-regensburg.de}
\date{To apprear in the special issue \\
``Extrapolation and Convergence Acceleration Methods'' \\
Journal of Computational and Applied Mathematics \\
(edited by Claude Brezinski)}

\maketitle

\begin{abstract} 
The prediction properties of Aitken's iterated $\Delta^2$ process,
Wynn's epsilon algorithm, and Brezinski's iterated theta algorithm for
(formal) power series are analyzed. As a first step, the defining
recursive schemes of these transformations are suitably rearranged in
order to permit the derivation of accuracy-through-order
relationships. On the basis of these relationships, the rational
approximants can be rewritten as a partial sum plus an appropriate
transformation term. A Taylor expansion of such a transformation term,
which is a rational function and which can be computed recursively,
produces the predictions for those coefficients of the (formal) power
series which were not used for the computation of the corresponding
rational approximant.
\end{abstract}

\setcounter{equation}{0}

\section{Introduction}
\label{Sec:Intro}

In applied mathematics and in theoretical physics, Pad\'{e} approximants
are now used almost routinely to overcome problems with slowly
convergent or divergent power series. Of course, there is an extensive
literature on Pad\'{e} approximants: In addition to countless articles,
there are several textbooks
\cite{Ba75,Ba90,BaGM96,Bre80,Bul97,CuWu87,DraVaIn87,Gil78,Po94}, review
articles \cite{Ba65,Ba72,Bas72,BrIs94,BrIs95,Gra72,ZiJ71}, collections 
of articles and proceedings \cite{BaGa70,Br91c,Ca76,Cu88,Cu94,dBrVRo81,%
GilPinSie85,GM73a,GM73b,GMSaVa84,SaVa77,WerBue84,Wuy79a}, bibliographies
\cite{Br76,Br91b,Wuy79b}, and there is even a book \cite{Br91a} and an 
article \cite{Bre96}, respectively, treating the history of Pad\'{e}
approximants and related topics. A long but by no means complete list of
applications of Pad\'{e} approximants in physics and chemistry can be 
found in Section 4 of \cite{We94b}.

The revival of the interest in Pad\'{e} approximants was initiated by
two articles by Shanks \cite{Sha55} and Wynn \cite{Wy56a},
respectively. These articles, which stimulated an enormous amount of
research, were published in 1956 at a time when electronic computers
started to become more widely available. Shanks \cite{Sha55} introduced
a sequence transformation which produces Pad\'{e} approximants if the
input data are the partial sums of a power series, and Wynn \cite{Wy56a}
showed that this transformation can be computed conveniently and
effectively by a recursive scheme now commonly called the epsilon
algorithm. As a consequence of the intense research initiated by Shanks
\cite{Sha55} and Wynn \cite{Wy56a}, the mathematical properties of
Pad\'{e} approximants are now fairly well understood, and it is
generally accepted that Pad\'{e} approximants are extremely useful
numerical tools which can be applied profitably in a large variety of
circumstances.

This intense research of course also showed that Pad\'{e} approximants
have certain limitations and shortcomings. For example, Pad\'{e}
approximants are in principle limited to convergent and divergent power
series and cannot help in the case of many other slowly convergent
sequences and series with different convergence types.

The convergence type of numerous practically important sequences $\{ s_n
\}_{n=0}^{\infty}$ can be classified by the asymptotic condition
\begin{equation}
\lim_{n \to \infty} \frac {s_{n+1} - s} {s_n - s} \; = \; \rho \, ,
\label{LimRatSeq}
\end{equation}
which closely resembles the well known ratio test for infinite
series. Here, $s = s_{\infty}$ is the limit of $\{ s_n
\}_{n=0}^{\infty}$ as $n \to \infty$. A convergent sequence satisfying
(\ref{LimRatSeq}) with $\vert \rho \vert < 1$ is called {\em linearly}
convergent, and it is called {\em logarithmically} convergent if $\rho =
1$. The partial sums of a power series with a nonzero, but finite radius
of convergence are a typical example of a linearly convergent sequence.
The partial sums of the Dirichlet series for the Riemann zeta
function,  
\begin{equation}
\zeta (z) \; = \; \sum_{m=0}^{\infty} \, (m+1)^{-z} \, , 
\qquad \mathrm{Re} (z) > 1 \, ,
\end{equation}
which is notorious for its extremely slow convergence if $\mathrm{Re}
(z)$ is only slightly larger than one, are a typical example of a
logarithmically convergent sequence.

Pad\'{e} approximants as well as the closely related epsilon algorithm
\cite{Wy56a} are known to accelerate effectively the convergence of 
linearly convergent power series and they are also able to sum many
divergent power series. However, they fail completely in the case of
logarithmic convergence (compare for example \cite[Theorem
12]{Wy66a}). Moreover, in the case of divergent power series whose
series coefficients grow more strongly than factorially, Pad\'{e}
approximants either converge too slowly to be numerically useful
\cite{CizVrs82,Sim82} or are not at all able to accomplish a summation
to a unique finite generalized limit \cite{GraGre78}. Consequently, the
articles by Shanks \cite{Sha55} and Wynn \cite{Wy56a} also stimulated
research on sequence transformations. The rapid progress in this field
is convincingly demonstrated by the large number of monographs and
review articles on sequence transformations which appeared in recent
years
\cite{Bre77,Bre78,Bre97,BreRZa91,De88,LieLueSh95,MarSha83,Wal96,We89,Wi81}.

In some, but by no means in all cases, sequence transformations are able
to do better than Pad\'{e} approximants, and it may even happen that
they clearly outperform Pad\'{e} approximants. Thus, it may well be
worth while to investigate whether it is possible to use instead of
Pad\'{e} approximants more specialized sequence transformations which
may be better adapted to the problem under consideration. For example,
the present author used sequence transformations successfully as
computational tools in such diverse fields as the evaluation of special
functions \cite{We89,We90,We94a,We94b,We96c,WeCi90,JeMoSoWe99,SteWe90},
the evaluation of molecular multicenter integrals of exponentially
decaying functions \cite{GroWeSte86,HoWe95,We94b,WeGroSte86,WeSte89},
the summation of strongly divergent quantum mechanical perturbation
expansions
\cite{CiViWe91,CiViWe93,CiWeBrSp96,We90,We92,We94b,WeCiVi91,WeCi90,%
We96a,We96b,We96d,We97}, and the extrapolation of quantum chemical
\emph{ab initio} calculations for oligomers to the infinite chain limit 
of quasi-onedimensional stereoregular polymers
\cite{We94b,WeLie90,CioWe93}. In vast majority of these applications, it
was either not possible to use Pad\'{e} approximants at all, or
alternative sequence transformations did a better job.

In most practical applications of Pad\'{e} approximants or also of
sequence transformations, the partial sums of (formal) power series are
transformed into rational approximants with the intention of either
accelerating convergence or to accomplish a summation to a finite
(generalized) limit in the case of divergence. Pad\'{e} approximants and
sequence transformations are normally not used for the computation of
the coefficients of the power series. In the majority of applications,
the computation of the coefficients of power series is not the most
serious computational problem, and conventional methods for the
computation of the coefficients usually suffice.

However, in the case of certain perturbation expansions as they for
instance occur in high energy physics, in quantum field theory, or in
quantum chromodynamics, the computational problems can be much more
severe. Not only do these perturbation expansions, which are power
series in some coupling constant, diverge quite strongly for every
nonzero value of the coupling constant, but it is also extremely
difficult to compute more than just a few of the perturbation series
coefficients. Moreover, due to the the complexity of the computations
and the necessity of making often drastic approximations, the
perturbation series coefficients obtained in this way are usually
affected by comparatively large relative errors. Under such adverse
circumstances, it has recently become customary to use Pad\'{e}
approximants to make predictions about the leading unknown coefficients
of perturbation expansions as well as to make consistency checks for the 
previously calculated coefficients
\cite{BroEllGarKarSam98,ChiEliSte99a,ChiEliSte99b,EllGarKarSam96a,%
EliSteChiMigSpr98,EllGarKarSam96b,EllJacJonKarSam98,EllKarSam97,Kar98,%
SamAbrYu97,SamEllKar95,SamLi94,SamLiSte93,SamLiSte95,SteEli98}. 

On a heuristic level, the prediction capability of Pad\'{e}
approximants, which was apparently first used by Gilewicz \cite{Gil73},
can be explained quite easily. Let us assume that a function $f$
possesses the following (formal) power series,
\begin{equation}
f (z) \; = \; \sum_{\nu=0}^{\infty} \, \gamma_{\nu} \, z^{\nu} \, ,
\label{PowSer_f}
\end{equation}
and that we want to transform the sequence of its partial sums
\begin{equation}
f_{n} (z) \; = \; \sum_{\nu=0}^{n} \, \gamma_{\nu} \, z^{\nu}
\label{ParSum_f}
\end{equation}
into a doubly indexed sequence of Pad\'{e} approximants 
\begin{equation}
[l/m]_f (z) \; = \; P_l (z)/Q_m (z) \, .
\label{Def_PA}
\end{equation}
As is well known \cite{Ba75,BaGM96}, the coefficients of the polynomials
$P_l (z) = p_0 + p_1 z + \ldots + p_l z^l$ and $Q_m (z) = 1 + q_1 z +
\ldots + q_m z^m$ are chosen in such a way that the Taylor expansion of
the Pad\'{e} approximant agrees as far as possible with the (formal)
power series (\ref{PowSer_f}):
\begin{equation}
f (z) - P_l (z)/Q_m (z) \; = \; O \bigl( z^{l+m+1} \bigr) \, , 
\qquad z \to 0 \, .
\label{PadeOrdEst}
\end{equation}
This \emph{accuracy-through-order} relationship implies that the
Pad\'{e} approximant to $f (z)$ can be written as the partial sum, from
which it was constructed, plus a term which was generated by the
transformation of the partial sum to the rational approximant:
\begin{equation}
[l/m]_f (z) \; = \; \sum_{\nu=0}^{l+m} \, \gamma_{\nu} \, z^{\nu} 
\, + \, z^{l+m+1} \, \mathcal{P}_{l}^{m} (z)
\; = \; f_{l+m} (z) \, + \, z^{l+m+1} \, \mathcal{P}_{l}^{m} (z) \, .
\end{equation}
Similarly, the (formal) power series (\ref{PowSer_f}) can be expressed
as follows:
\begin{equation}
f (z) \; = \; \sum_{\nu=0}^{l+m} \, \gamma_{\nu} \, z^{\nu} 
\, + \, z^{l+m+1} \, \mathcal{F}_{l+m+1} (z) \; = \; 
f_{l+m} (z) \, + \, z^{l+m+1} \, \mathcal{F}_{l+m+1} (z) \, .
\end{equation}

Let us now assume that the Pad\'{e} approximant $[l/m]_f (z)$ provides a
sufficiently accurate approximation to $f (z)$. Then, the Pad\'{e}
transformation term $\mathcal{P}_{l}^{m} (z)$ must also provide a
sufficiently accurate approximation to the truncation error
$\mathcal{F}_{l+m+1} (z)$ of the (formal) power series. In general, we
have no reason to assume that $\mathcal{P}_{l}^{m} (z)$ could be equal
to $\mathcal{F}_{l+m+1} (z)$ for finite values of $l$ and
$m$. Consequently, Taylor expansions of $\mathcal{P}_{l}^{m} (z)$ and
$\mathcal{F}_{l+m+1} (z)$, respectively, will in general produce
different results. Nevertheless, the \emph{leading} coefficients of the
Taylor expansion for $\mathcal{P}_{l}^{m} (z)$ should provide
sufficiently accurate approximations to the corresponding coefficients
of the Taylor series for $\mathcal{F}_{l+m+1} (z)$.

It is important to note that this prediction capability does not depend
on the convergence of the power series expansions for
$\mathcal{P}_{l}^{m} (z)$ and $\mathcal{F}_{l+m+1} (z)$,
respectively. Pad\'{e} approximants are able to make predictions about
series coefficients even if the power series (\ref{PowSer_f}) for $f$ as
well as the power series expansions for $\mathcal{P}_{l}^{m}$ and
$\mathcal{F}_{l+m+1} (z)$ are only asymptotic as $z \to 0$. This fact
explains why the prediction capability of Pad\'{e} approximants can be
so very useful in the case of violently divergent perturbation
expansions.

Let us now assume that a sequence transformation also produces a
convergent sequence of rational approximants if it acts on the partial
sums (\ref{ParSum_f}) of the (formal) power series
(\ref{PowSer_f}). Then, by the same line of reasoning, these rational
approximants should also be able to make predictions about the leading
coefficients of the power series, which were not used for the
construction of the rational approximant. It seems that these ideas were
first formulated by Sidi and Levin \cite{SidLev83} and Brezinski
\cite{Bre85}. Recently, these ideas were extended by Pr\'{e}vost and
Vekemans \cite{PreVeke99} who discussed prediction methods for sequences
which they called $\varepsilon_p$ and partial Pad\'{e} prediction,
respectively. Moreover, in \cite{We97} it was shown that suitably chosen
sequence transformations can indeed make more accurate predictions about
unknown power series coefficients than Pad\'e approximants.

Consequently, it should be interesting to analyze the prediction
properties of sequence transformations. In this this article, only
Aitken's iterated $\Delta^2$ algorithm, Wynn's epsilon algorithm and the
iteration of Brezinski's theta algorithm will be considered.  Further
studies on the prediction properties of other sequence transformations
are in progress and will be presented elsewhere.

If the prediction properties of sequence transformations are to be
studied, there is an additional complication which is absent in the case
of Pad\'{e} approximants. The accuracy-through-order relationship
(\ref{PadeOrdEst}) leads to a system of $l+m+1$ linear equations for the
coefficients of the polynomials $P_l (z) = p_0 + p_1 z + \ldots + p_l
z^l$ and $Q_m (z) = 1 + q_1 z + \ldots + q_m z^m$ of the Pad\'{e}
approximant (\ref{Def_PA}) \cite{Ba75,BaGM96}. If this system of
equations has a solution, then it is automatically guaranteed that the
Pad\'{e} approximant obtained in this way satisfies the
accuracy-through-order relationship (\ref{PadeOrdEst}).

In the case of the sequence transformations considered in this article,
the situation is in general more complicated. These transformations are
not defined as solutions of systems of linear equations, but via
nonlinear recursive schemes. Moreover, their accuracy-through-order
relationships are with the exception of Wynn's epsilon algorithm unknown
and have to be derived via their defining recursive schemes.

On the basis of these accuracy-through-order relationships, it is
possible to construct explicit recursive schemes for the transformation
errors as well as for the first coefficient of the power series which
was not used for the computation of the rational approximant.

In Section \ref{Sec:Aitken}, the the accuracy-through-order and
prediction properties of Aitken's iterated $\Delta^2$ process are
analyzed. In Section \ref{Sec:EpsAl}, the analogous properties of Wynn's
epsilon algorithm are discussed, and in Section \ref{Sec:ThetIt},
Brezinski's iterated theta algorithm is treated. In Section
\ref{Sec:Applics}, some applications of the new results are
presented. This article is concluded by Section \ref{Sec:SumCon} which
contains a short summary.

\setcounter{equation}{0}

\section{Aitken's Iterated $\Delta^2$ Process}
\label{Sec:Aitken}

Let us consider the following model sequence:
\begin{equation}
s_n \; = \; s \, + \, c \, \lambda^n \, , \qquad c \ne 0, 
\quad \vert \lambda \vert \ne 1 \, , \quad n \in \mathbb{N}_0 \, .
\label{AitModSeq}
\end{equation}
For $n \to \infty$, this sequence obviously converges to its limit $s$
if $0 < \vert \lambda \vert < 1$, and it diverges away from its
generalized limit $s$ if $\vert \lambda \vert > 1$.

A sequence transformation, which is able to determine the (generalized)
limit $s$ of the model sequence (\ref{AitModSeq}) from the numerical
values of three consecutive sequence elements $s_n$, $s_{n+1}$ and
$s_{n+2}$, can be constructed quite easily. Just consider $s$, $c$, and
$\lambda$ as unknowns of the linear system $s_{n+j} = s + c
\lambda^{n+j}$ with $j = 0, 1, 2$. A short calculation shows that
\begin{equation}
\mathcal{A}_{1}^{(n)} \; = \; s_n \, - \, 
\frac{[\Delta s_n]^2}{\Delta^2 s_n} \, ,  \qquad n \in \mathbb{N}_0 \, ,
\label{AitFor_1}
\end{equation}
is able to determine the (generalized) limit of the model sequence
(\ref{AitModSeq}) according to $\mathcal{A}_{1}^{(n)} = s$. It should be
noted that $s$ can be determined in this way, no matter whether the
sequence (\ref{AitModSeq}) converges or diverges. The forward
difference operator $\Delta$ in (\ref{AitFor_1}) is defined by its
action on a function $g = g (n)$:
\begin{equation}
\Delta g (n) \; = \; g (n+1) \, - \, g (n) \, .
\end{equation}

The $\Delta^2$ formula (\ref{AitFor_1}) is certainly one of the oldest
sequence transformations. It is usually attributed to Aitken
\cite{Ai26}, but it is actually much older. Brezinski
\cite[pp.\ 90 - 91]{Br91a} mentioned that in 1674 Seki Kowa, the
probably most famous Japanese mathematician of that period, tried to
obtain better approximations to $\pi$ with the help of this $\Delta^2$
formula, and according to Todd \cite[p.\ 5]{Tod62} it was in principle
already known to Kummer \cite{Kum1837}.

There is an extensive literature on Aitken's $\Delta^2$ process. For
example, it was discussed by Lubkin \cite{Lub52}, Shanks \cite{Sha55},
Tucker \cite{Tuc67,Tuc69}, Clark, Gray, and Adams \cite{ClGrAd69},
Cordellier \cite{Cor79}, Jurkat \cite{Jur83}, Bell and Phillips
\cite{BelPhi84}, and Weniger \cite[Section 5]{We89}. A multidimensional
generalization of Aitken's transformation to vector sequences was
discussed by MacLeod \cite{MacL86}. Modifications and generalizations of
Aitken's $\Delta^2$ process were proposed by Drummond \cite{Dru76},
Jamieson and O'Beirne \cite{JamOBe78}, Bj{\o}rstad, Dahlquist, and
Grosse \cite{BjoDahGro81}, and Sablonniere \cite{Sab92}. Then, there is
a close connection between the Aitken process and Fibonacci numbers, as
discussed by McCabe and Phillips \cite{McCPhi85} and Arai, Okamoto, and
Kametaka \cite{ArOkKa88}. The properties of Aitken's $\Delta^2$ process
are also discussed in books by Baker and Graves-Morris \cite{BaGM96},
Brezinski \cite{Bre77,Bre78}, Brezinski and Redivo Zaglia
\cite{BreRZa91}, Delahaye \cite{De88}, Walz \cite{Wal96}, and 
Wimp \cite{Wi81}.

The power of Aitken's $\Delta^2$ process is of course limited since it
is designed to eliminate only a single exponential term. However, its
power can be increased considerably by iterating it, yielding the
following nonlinear recursive scheme:
\begin{subequations}
\label{ItAit_1}
\begin{eqnarray}
\mathcal{A}_{0}^{(n)} & = & s_n \, , \qquad n \in \mathbb{N}_0 \, ,
\\
\mathcal{A}_{k+1}^{(n)} & = & \mathcal{A}_{k}^{(n)} \, - \, 
\frac{\bigl[ \Delta \mathcal{A}_{k}^{(n)} \bigr]^2}
{\Delta^2 \mathcal{A}_{k}^{(n)}} \, , \qquad k, n \in \mathbb{N}_0 \, .
\end{eqnarray}
\end{subequations}
In the case of doubly indexed quantities like $\mathcal{A}_{k}^{(n)}$,
it will always be assumed that the difference operator $\Delta$ only
acts on the superscript $n$ but not on the subscript $k$:
\begin{equation}
\Delta \mathcal{A}_{k}^{(n)} \; = \; 
\mathcal{A}_{k}^{(n+1)} \, - \, \mathcal{A}_{k}^{(n)} \, .
\end{equation}

The numerical performance of Aitken's iterated $\Delta^2$ process was
studied in \cite{We89,SmiFor82}. Concerning the theoretical properties
of Aitken's iterated $\Delta^2$ process, very little seems to be
known. Hillion \cite{Hil75} was able to find a model sequence for which
the iterated $\Delta^2$ process is exact. He also derived a
determinantal representation for $\mathcal{A}_k^{(n)}$. However,
Hillion's expressions for $\mathcal{A}_k^{(n)}$ contain explicitly the
lower order transforms $\mathcal{A}_0^{(n)}, \ldots,
\mathcal{A}_{k-1}^{(n)}, \ldots, \mathcal{A}_0^{(n+k)}, \ldots, 
\mathcal{A}_{k-1}^{(n+k)}$. Consequently, 
it seems that Hillion's result \cite{Hil75} -- although interesting from
a formal point of view -- cannot help much to analyze the prediction
properties of $\mathcal{A}_{k}^{(n)}$.

If we want to use Aitken's iterated $\Delta^2$ process for the
prediction of unknown series coefficients, we first have to derive its
accuracy-through-order relationship of the type of (\ref{PadeOrdEst}) on
the basis of the recursive scheme (\ref{ItAit_1}).

It is a direct consequence of the recursive scheme (\ref{ItAit_1}) that
$2k+1$ sequence elements $s_{n}$, $s_{n+1}$, \ldots, $s_{n+2k}$ are
needed for the computation of $\mathcal{A}_{k}^{(n)}$. Thus, we now
choose as input data the partial sums (\ref{ParSum_f}) of the (formal)
power series (\ref{PowSer_f}) according to $s_n = f_n (z)$, and
conjecture that all coefficients $\gamma_{0}$, $\gamma_{1}$, \ldots,
$\gamma_{n+2k}$, which were used for the construction of
$\mathcal{A}_{k}^{(n)}$, are exactly reproduced by a Taylor
expansion. This means that we have to look for an accuracy-through-order
relationship of the following kind:
\begin{equation}
f (z) \, - \, \mathcal{A}_{k}^{(n)} \; = \; 
O \bigl( z^{n+2k+1} \bigr) \, , \qquad z \to 0 \, .
\label{Ait_OrdEst}
\end{equation}
Such an accuracy-through-order relationship would imply that
$\mathcal{A}_{k}^{(n)}$ can be expressed as follows:
\begin{equation}
\mathcal{A}_{k}^{(n)} \; = \; f_{n+2k} (z) \, + \, G_{k}^{(n)} \, 
z^{n+2k+1} \, + \, O \bigl( z^{n+2k+2} \bigr) \, , 
\qquad z \to 0 \, .
\label{AitAsyZero}
\end{equation}
The constant $G_{k}^{(n)}$ is the prediction made for the coefficient
$\gamma_{n+2k+1}$, which is the first coefficient of the power series
(\ref{PowSer_f}) not used for the computation of
$\mathcal{A}_{k}^{(n)}$. 
 
Unfortunately, the recursive scheme (\ref{ItAit_1}) is not suited for
our purposes. This can be shown by computing $\mathcal{A}_{1}^{(n)}$
from the partial sums $f_{n} (z)$, $f_{n+1} (z)$, and $f_{n+2} (z)$:
\begin{equation}
\mathcal{A}_{1}^{(n)} \; = \; f_{n} (z) \, + \, 
\frac{\bigl[ \gamma_{n+1} \bigr]^2 z^{n+1}}
{\gamma_{n+1} - \gamma_{n+2} z} \, .
\label{Ait_1_1}
\end{equation}
Superficially, it looks as if $\mathcal{A}_{1}^{(n)}$ is not of the type
of (\ref{AitAsyZero}). However, the rational expression on the
right-hand side contains the missing terms $\gamma_{n+1} z^{n+1}$ and
$\gamma_{n+2} z^{n+2}$. We only have to use $1/(1-y) = 1 + y +
y^2/(1-y)$ with $y = \gamma_{n+2} z / \gamma_{n+1}$ to obtain an
equivalent expression with the desired features:
\begin{equation}
\mathcal{A}_{1}^{(n)} \; = \; f_{n+2} (z) \, + \, 
\frac{\bigl[ \gamma_{n+2} \bigr]^2 z^{n+3}}
{\gamma_{n+1} - \gamma_{n+2} z} \, .
\label{Ait_1_2}
\end{equation}
Thus, an expression, which is in agreement with (\ref{AitAsyZero}), can
be obtained easily in the case of the simplest transform
$\mathcal{A}_{1}^{(n)}$. Moreover, (\ref{Ait_1_2}) makes the prediction
$G_{1}^{(n)} = \bigl[ \gamma_{n+2} \bigr]^2/\gamma_{n+1}$ for the first
series coefficient $\gamma_{n+3}$ not used for the computation of
$\mathcal{A}_{1}^{(n)}$. Of course, by expanding the denominator on the
right-hand side of (\ref{Ait_1_2}) further predictions on series
coefficients with higher indices can be made.

In the case of more complicated transforms $\mathcal{A}_{k}^{(n)}$ with
$k > 1$, it is by no means obvious whether and how the necessary
manipulations, which would transform an expression of the type of
(\ref{Ait_1_1}) into an expression of the type of (\ref{Ait_1_2}), can
be done. Consequently, it is advantageous to replace the recursive
scheme (\ref{ItAit_1}) by an alternative recursive scheme, which
directly leads to appropriate expressions for $\mathcal{A}_{k}^{(n)}$
with $k > 1$.

Many different expressions for $\mathcal{A}_{1}^{(n)}$ in terms of
$s_n$, $s_{n+1}$, and $s_{n+2}$ are known \cite[Section 5.1]{We89}.
These expressions are all mathematically equivalent although their
numerical properties may differ. Comparison with (\ref{Ait_1_2}) shows
that the for our purposes appropriate expression is
\cite[Eq.\ (5.1-7)]{We89}
\begin{equation}
\mathcal{A}_{1}^{(n)} \; = \;
s_{n+2} \, - \, \frac{[\Delta s_{n+1}]^2}{\Delta^2 s_n} \, .
\label{AitFor_2}
\end{equation}
Just like (\ref{AitFor_1}), this expression can be iterated and yields
\begin{subequations}
\label{ItAit_2}
\begin{eqnarray}
\mathcal{A}_{0}^{(n)} & = & s_n \, , \qquad n \in \mathbb{N}_0 \, ,
\label{ItAit_2_a}
\\
\mathcal{A}_{k+1}^{(n)} & = & \mathcal{A}_{k}^{(n+2)} \, - \, 
\frac{\bigl[ \Delta \mathcal{A}_{k}^{(n+1)} \bigr]^2}
{\Delta^2 \mathcal{A}_{k}^{(n)}} \, , \qquad k, n \in \mathbb{N}_0 \, .
\label{ItAit_2_b}
\end{eqnarray}
\end{subequations}
The recursive schemes (\ref{ItAit_1}) and (\ref{ItAit_2}) are
mathematically completely equivalent. However, for our purposes -- the
analysis of the prediction properties of Aitken's iterated $\Delta^2$
process in the case of power series -- the recursive scheme
(\ref{ItAit_2}) is much better suited.

Next, we rewrite the partial sums (\ref{ParSum_f}) of the (formal)
power series (\ref{PowSer_f}) according to
\begin{equation}
f_n (z) \; = \; f (z) \, - \, 
\sum_{\nu=0}^{\infty} \, \gamma_{n+\nu+1} \, z^{n+\nu+1}
\label{ParSum_Rem}
\end{equation}
and use them as input data in the recursive scheme (\ref{ItAit_2}). This
yields the following expression:
\begin{equation}
\mathcal{A}_{k}^{(n)} \; = \; 
f (z) \, + \, z^{n+2k+1} \, R_{k}^{(n)} (z) \, ,
\qquad k, n \in \mathbb{N}_0 \, .	
\label{Ait_AccThrOrd}
\end{equation}
The quantities $R_{k}^{(n)} (z)$ can be computed with the help of the
following recursive scheme which is a direct consequence of the
recursive scheme (\ref{ItAit_2}) for $\mathcal{A}_{k}^{(n)}$:
\begin{subequations}
\label{RecR}
\begin{eqnarray}
R_{0}^{(n)} (z) & = & 
- \, \sum_{\nu=0}^{\infty} \, \gamma_{n+\nu+1} \, z^{\nu} \; = \;
\frac{f_{n} (z) - f (z)}{z^{n+1}} \, , \qquad n \in \mathbb{N}_0 \, ,
\\
R_{k+1}^{(n)} (z) & = & R_{k}^{(n+2)} (z) \, - \, 
\frac{\bigl[ \delta R_{k}^{(n+1)} (z) \bigr]^2}
{\delta^2 R_{k}^{(n)} (z)} \, , \qquad k, n \in \mathbb{N}_0 \, .
\end{eqnarray}
\end{subequations}
In (\ref{RecR}), we use the shorthand notation
\begin{subequations}
\label{delta_X}
\begin{eqnarray}
\delta X_{k}^{(n)} (z) & = & 
z X_{k}^{(n+1)} (z) \, - \, X_{k}^{(n)} (z) \, , 
\\
\delta^2 X_{k}^{(n)} (z) & = & 
z \delta X_{k}^{(n+1)} (z) \, - \, \delta X_{k}^{(n)} (z) 
\nonumber \\
& = & z^2 X_{k}^{(n+2)} (z) \, - \, 2 z X_{k}^{(n+1)} (z) \, + \, 
X_{k}^{(n)} (z) \, .
\end{eqnarray}
\end{subequations}

It seems that we have now accomplished our aim since
(\ref{Ait_AccThrOrd}) has the right structure to serve as an
accuracy-through-order relationship for Aitken's iterated $\Delta^2$
process. Unfortunately, this conclusion is in general premature and we
have to require that the input data satisfy some additional
conditions. One must not forget that Aitken's $\Delta^2$ formula
(\ref{AitFor_2}) as well as its iteration (\ref{ItAit_2}) cannot be
applied to arbitrary input data. One obvious potential complication,
which has to be excluded, is that (\ref{ItAit_2_b}) becomes undefined if
$\Delta^2 \mathcal{A}_{k}^{(n)} = 0$. Thus, if we want to transform the
partial sums (\ref{ParSum_f}) of the (formal) power series
(\ref{PowSer_f}), it is natural to require that all series coefficients
are nonzero, i.e., $\gamma_{\nu} \ne 0$ for all $\nu \in \mathbb{N}_0$.

Unfortunately, this is only a minimal requirement and not yet enough for
our purposes. If $z^{n+2k+1} R_{k}^{(n)} (z)$ in (\ref{Ait_AccThrOrd})
is to be of order $O \bigl( z^{n+2k+1} \bigr)$ as $z \to 0$, then the
$z$-independent part $C_{k}^{(n)}$ of $R_{k}^{(n)} (z)$ defined by
\begin{equation}
R_{k}^{(n)} (z) \; = \; C_{k}^{(n)} \, + \, O (z) \, , 
\qquad z \to 0 \, ,
\label{Ait_z_indep_NZ_TruncErr_a}
\end{equation} 
has to satisfy
\begin{equation}
C_{k}^{(n)} \; \ne \; 0 \, , \qquad k, n \in \mathbb{N}_0 \, .
\label{Ait_z_indep_NZ_TruncErr_b}
\end{equation}
If these conditions are satisfied, we can be sure that
(\ref{Ait_AccThrOrd}) is indeed the accuracy-through-order relationship
we have been looking for.

Personally, I am quite sceptical that it would be easy to characterize
\emph{theoretically} those power series which give rise to truncation 
errors $R_{k}^{(n)} (z)$ satisfying (\ref{Ait_z_indep_NZ_TruncErr_a}) and
(\ref{Ait_z_indep_NZ_TruncErr_b}). Fortunately, it can easily be checked
\emph{numerically} whether a given (formal) power series leads to 
truncation errors whose $z$-independent parts are nonzero. If we set $z
= 0$ in (\ref{RecR}) and use (\ref{Ait_z_indep_NZ_TruncErr_a}), we obtain
the following recursive scheme:
\begin{subequations}
\label{RecC}
\begin{eqnarray}
C_{0}^{(n)} & = & - \, \gamma_{n+1} \, , \qquad n \in \mathbb{N}_0 \, ,
\\
C_{k+1}^{(n)} & = & C_{k}^{(n+2)} \, - \, 
\frac{\bigl[ C_{k}^{(n+1)} \bigr]^2}{C_{k}^{(n)}} \, ,
\qquad k, n \in \mathbb{N}_0 \, .
\end{eqnarray}
\end{subequations}

Let us now assume that we know for a given (formal) power series that
the $z$-independent parts $C_{k}^{(n)}$ of the truncation errors
$R_{k}^{(n)} (z)$ in (\ref{Ait_AccThrOrd}) are nonzero -- either from a
mathematical proof or from a brute force calculation using
(\ref{RecC}). Then, (\ref{Ait_AccThrOrd}) is indeed the
accuracy-through-order relationship we have been looking for, which
implies that $\mathcal{A}_{k}^{(n)}$ can be expressed as follows:
\begin{equation}
\mathcal{A}_{k}^{(n)} \; = \; f_{n+2k} (z) \, + \, 
z^{n+2k+1} \, \Phi_{k}^{(n)} (z) \, , \qquad k, n \in \mathbb{N}_0 \, .
\label{AitRemPhi}
\end{equation}
If we use this ansatz in (\ref{ItAit_2}), we obtain the following
recursive scheme:
\begin{subequations}
\label{RecPhi}
\begin{eqnarray}
\Phi_{0}^{(n)} (z) & = & 0 \, , \qquad n \in \mathbb{N}_0 \, ,
\label{RecPhi_a}
\\
\Phi_{k+1}^{(n)} (z) & = & \Phi_{k}^{(n+2)} (z) \, - \, 
\frac{\bigl[ \gamma_{n+2k+2} + \delta \Phi_{k}^{(n+1)} (z) \bigr]^2}
{\gamma_{n+2k+2} z - \gamma_{n+2k+1} + \delta^2 \Phi_{k}^{(n)} (z)} \, ,
\qquad k, n \in \mathbb{N}_0 \, .
\label{RecPhi_b}
\end{eqnarray}
\end{subequations}
Here, $\delta \Phi_{k}^{(n)} (z)$ and $\delta^2 \Phi_{k}^{(n)} (z)$ are
defined by (\ref{delta_X}). For $k = 0$, (\ref{RecPhi_b}) yields
\begin{equation}
\Phi_{1}^{(n)} (z) \; = \; \frac{\bigl[ \gamma_{n+2} \bigr]^2}
{\gamma_{n+1} - \gamma_{n+2} z} \, ,
\end{equation}
which is in agreement with (\ref{Ait_1_2}).

A comparison of (\ref{AitAsyZero}) and (\ref{AitRemPhi}) yields
\begin{equation}
\Phi_{k}^{(n)} (z) \; = \; G_{k}^{(n)} \, + \, O \bigl( z \bigr) \, ,
\qquad z \to 0 \, .
\label{Phi2G}
\end{equation}
Consequently, the $z$-independent part $G_{k}^{(n)}$ of $\Phi_{k}^{(n)}
(z)$ is the prediction for the first coefficient $\gamma_{n+2k+1}$ not
used for the computation of $\mathcal{A}_{k}^{(n)}$.

If we set $z = 0$ in the recursive scheme (\ref{RecPhi}) and use
(\ref{Phi2G}), we obtain the following recursive scheme for the
predictions $G_{k}^{(n)}$:
\begin{subequations}
\label{RecG}
\begin{eqnarray}
G_{0}^{(n)} & = & 0 \, , \qquad n \in \mathbb{N}_0 \, ,
\\
G_{1}^{(n)} & = & \bigl[ \gamma_{n+2} \bigr]^2/\gamma_{n+1} \, ,
\qquad n \in \mathbb{N}_0 \, ,
\\
G_{k+1}^{(n)} & = & G_{k}^{(n+2)} \, + \, \frac
{\bigl[ \gamma_{n+2k+2} - G_{k}^{(n+1)} \bigr]^2}
{\gamma_{n+2k+1} - G_{k}^{(n)}} \, , \qquad k, n \in \mathbb{N}_0 \, .
\end{eqnarray}
\end{subequations}

The $z$-independent parts $C_{k}^{(n)}$ of $R_{k}^{(n)} (z)$ and
$G_{k}^{(n)}$ of $\Phi_{k}^{(n)} (z)$, respectively, are connected.  A
comparison of (\ref{Ait_AccThrOrd}), (\ref{Ait_z_indep_NZ_TruncErr_a}),
(\ref{AitRemPhi}), and (\ref{Phi2G}) yields:
\begin{equation}
G_{k}^{(n)} \; = \; C_{k}^{(n)} \, + \, \gamma_{n+2k+1} \, .
\end{equation}

In this article, rational approximants will always be used in such a way
that the input data -- the partial sums (\ref{ParSum_f}) of the (formal)
power series (\ref{PowSer_f}) -- are computed in an outer loop, and for
each new partial sum a new approximation to the limit is calculated. If
the index $m$ of the last partial sum $f_{m} (z)$ is even, $m = 2 \mu$,
we use in the case of Aitken's iterated $\Delta^2$ process as
approximation to the limit $f (z)$ the transformation
\begin{equation}
\bigl\{ f_{0} (z), f_{1} (z), \ldots, f_{2\mu} (z) \bigr\} \mapsto 
\mathcal{A}_{\mu}^{(0)} \, ,
\end{equation}
and if $m$ is odd, $m = 2\mu+1$, we use the transformation
\begin{equation}
\bigl\{ f_{1} (z), f_{2} (z), \ldots, f_{2\mu+1} (z) \bigr\} \mapsto 
\mathcal{A}_{\mu}^{(1)} \, .
\end{equation}
With the help of the notation $\Ent{x}$ for the integral part of $x$,
which is the largest integer $\nu$ satisfying the inequality $\nu \le
x$, these two relationships can be combined into a single equation,
yielding \cite[Eq.\ (5.2-6)]{We89}
\begin{equation}
\bigl\{ f_{m-2\Ent{m/2}} (z), f_{m-2\Ent{m/2}+1} (z), \ldots, 
f_{m} (z) \bigr\} \mapsto \mathcal{A}_{\Ent{m/2}}^{(m-2\Ent{m/2})} \, , 
\qquad m \in \mathbb{N}_0 \, .
\label{AitAppLim}
\end{equation}
The same strategy will also be used if for example the rational
expressions $R_{k}^{(n)} (z)$ defined by (\ref{Ait_AccThrOrd}) are
listed in a Table. This means that the $R_{k}^{(n)} (z)$ will also be
listed according to (\ref{AitAppLim}). The only difference is that the
$R_{k}^{(n)} (z)$ use as input data not the partial sums $f_{n} (z)$ but
the remainders $[f_{n} (z) - f (z)]/z^{n+1}$.

\setcounter{equation}{0}

\section{Wynn's Epsilon Algorithm}
\label{Sec:EpsAl}

Wynn's epsilon algorithm \cite{Wy56a} is the following nonlinear
recursive scheme:
\begin{subequations}
\label{eps_al}
\begin{eqnarray}
\epsilon_{-1}^{(n)} & \; = \; & 0 \, ,  
\qquad \epsilon_0^{(n)} \, = \, s_n \, , 
\qquad  n \in \mathbb{N}_0 \, , \\
\epsilon_{k+1}^{(n)} & \; = \; & \epsilon_{k-1}^{(n+1)} \, + \,
1 / [\epsilon_{k}^{(n+1)} - \epsilon_{k}^{(n)} ] \, ,
\qquad k, n \in \mathbb{N}_0 \, .
\end{eqnarray}
\end{subequations}  
The elements $\epsilon_{2k}^{(n)}$ with \emph{even} subscripts provide
approximations to the (generalized) limit $s$ of the sequence $\{ s_n
\}_{n=0}^{\infty}$ to be transformed, whereas the elements
$\epsilon_{2k+1}^{(n)}$ with \emph{odd} subscripts are only auxiliary
quantities which diverge if the whole process converges.

If the input data are the partial sums (\ref{ParSum_f}) of the (formal)
power series (\ref{PowSer_f}), $s_n = f_n (z)$, then Wynn \cite{Wy56a}
could show that his epsilon algorithm produces Pad\'{e} approximants:
\begin{equation}
\epsilon_{2 k}^{(n)} \; = \; [ n + k / k ]_f (z) \, .
\label{Eps_Pade}
\end{equation}

The epsilon algorithm is a close relative of Aitken's iterated
$\Delta^2$ process, and they have similar properties in convergence
acceleration and summation processes. A straightforward calculation
shows that $\mathcal{A}_{1}^{(n)} = \epsilon_{2}^{(n)}$. Hence, Aitken's
iterated $\Delta^2$ process may also be viewed as an iteration of
$\epsilon_{2}^{(n)}$. However, for $k > 1$, $\mathcal{A}_{k}^{(n)}$ and
$\epsilon_{2k}^{(n)}$ are in general different.

There is an extensive literature on the epsilon algorithm.  On p.\ 120
of Wimps book \cite{Wi81} it is mentioned that over 50 articles on the
epsilon algorithm were published by Wynn alone, and at least 30 articles
by Brezinski. As a fairly complete source of references Wimp recommends
Brezinski's first book \cite{Bre77}. However, this book was published in
1977, and since then many more articles on the epsilon algorithm have
been published. Consequently, any attempt to produce something
resembling a reasonably complete bibliography of Wynn's epsilon
algorithm would clearly be beyond the scope of this article.

In spite of its numerous advantageous features, Wynn's epsilon algorithm
(\ref{eps_al}) is not suited for our purposes. If the input data are
the partial sums (\ref{ParSum_f}) of the (formal) power series
(\ref{PowSer_f}), the accuracy-through-order relationship
(\ref{PadeOrdEst}) of Pad\'{e} approximants in combination with
(\ref{Eps_Pade}) implies that the elements of the epsilon table with
even subscripts can be expressed as
\begin{equation}
\epsilon_{2k}^{(n)} \; = \; f_{n+2k} (z) \, + \, g_{2k}^{(n)} \, 
z^{n+2k+1} \, + \, O \bigl( z^{n+2k+2} \bigr) \, , 
\qquad z \to 0 \, .
\label{EpsAsyZero}
\end{equation}
The constant $g_{2k}^{(n)}$ is the prediction made for the coefficient
$\gamma_{n+2k+1}$, which is the first coefficient of the power series
(\ref{PowSer_f}) not used for the computation of $\epsilon_{2k}^{(n)}$.

If we compute $\epsilon_{2}^{(n)}$ from the partial sums $f_{n} (z)$,
$f_{n+1} (z)$, and $f_{n+2} (z)$, we obtain because of
$\mathcal{A}_{1}^{(n)} = \epsilon_{2}^{(n)}$ the same expressions as in
the last section. Thus, we obtain a result which does not seem to be in
agreement with the accuracy-through-order relationship
(\ref{EpsAsyZero}):
\begin{equation}
\epsilon_{2}^{(n)} \; = \; f_{n+1} (z) \, + \, 
\frac{\gamma_{n+1} \gamma_{n+2} z^{n+2}}
{\gamma_{n+1} - \gamma_{n+2} z} \, .
\label{eps_2_1}
\end{equation}
Of course, the missing term $\gamma_{n+2} z^{n+2}$ can easily be
extracted from the rational expression on the right-hand side.  We only
have to use $1/(1-y) = 1 + y/(1-y)$ with $y = \gamma_{n+2} z /
\gamma_{n+1}$ to obtain as in the case of Aitken's iterated $\Delta^2$
algorithm an expression with the desired features:
\begin{equation}
\epsilon_{2}^{(n)} \; = \; f_{n+2} (z) \, + \, 
\frac{\bigl[ \gamma_{n+2} \bigr]^2 z^{n+3}}
{\gamma_{n+1} - \gamma_{n+2} z} \, .
\label{eps_2_2}
\end{equation}

This example shows that the accuracy-through-order relationship
(\ref{PadeOrdEst}) of Pad\'{e} approximants is by no means immediately
obvious from the epsilon algorithm (\ref{eps_al}). A further
complication is that the epsilon algorithm involves the elements
$\epsilon_{2k+1}^{(n)}$ with odd subscripts. These are only auxiliary
quantities which diverge if the whole process converges. Nevertheless,
they make it difficult to obtain order estimates and to reformulate the
epsilon algorithm in such a way that it automatically produces suitable
expressions for $\epsilon_{2k}^{(n)}$ of the type of (\ref{eps_2_2}).

The starting point for the construction of an alternative recursive
scheme, which would be suited for our purposes, is Wynn's cross rule
\cite[Eq.\ (13)]{Wy66b}:
\begin{eqnarray}
\lefteqn
{\left\{ \epsilon_{2k+2}^{(n)} - \epsilon_{2k}^{(n+1)} \right\}^{-1} 
\, + \,
\left\{ \epsilon_{2k-2}^{(n+2)} - \epsilon_{2k}^{(n+1)} \right\}^{-1}} 
\nonumber \\
& \qquad  = &
\left\{ \epsilon_{2k}^{(n)} - \epsilon_{2k}^{(n+1)} \right\}^{-1} 
\, + \,
\left\{ \epsilon_{2k}^{(n+2)} - \epsilon_{2k}^{(n+1)} \right\}^{-1} \, .
\label{CrossRule}
\end{eqnarray}
This expression permits the recursive computation of the elements
$\epsilon_{2k}^{(n)}$ with even subscripts without having to compute the
auxiliary quantities $\epsilon_{2k+1}^{(n)}$ with odd subscripts. The
price, one has to pay, is that the cross rule (\ref{CrossRule}) has a
more complicated structure than the extremely simple epsilon algorithm
(\ref{eps_al}).

A further complication is that for $k = 0$ the undefined element
$\epsilon_{-2}^{(n)}$ occurs in (\ref{CrossRule}). However, we obtain
results that are consistent with Wynn's epsilon algorithm (\ref{eps_al})
if we set $\epsilon_{-2}^{(n)} = \infty$.

Hence, instead of the epsilon algorithm (\ref{eps_al}), we can also use
the following recursive scheme:
\begin{subequations}
\label{CrossRule_1}
\begin{eqnarray}
\epsilon_{-2}^{(n)} & = & \infty \, , 
\qquad \epsilon_{0}^{(n)} \; = \; s_n \, , \qquad n \in \mathbb{N}_0 \, ,
\\
\epsilon_{2k+2}^{(n)} & = & \epsilon_{2k}^{(n+1)} \, + \,
\frac{\displaystyle 1}{\displaystyle 
\frac{1 \rule{0pt}{11pt}}{\displaystyle \Delta \epsilon_{2k}^{(n+1)}} - 
\frac{1}{\displaystyle \Delta \epsilon_{2k}^{(n)}} +
\frac{1}{\displaystyle 
\epsilon_{2k}^{(n+1)} - \epsilon_{2k-2}^{(n+2)}}} \, , 
\qquad k, n \in \mathbb{N}_0 \, .
\end{eqnarray}
\end{subequations}

For our purposes, this recursive scheme is an improvement over the
epsilon algorithm (\ref{eps_al}) since it does not contain the elements
$\epsilon_{2k+1}^{(n)}$ with odd subscripts. Nevertheless, it is not yet
what we need. The use of (\ref{CrossRule_1}) for the computation of
$\epsilon_{2}^{(n)}$ would produce (\ref{eps_2_1}) but not
(\ref{eps_2_2}). Fortunately, (\ref{CrossRule_1}) can easily be modified
to yield a recursive scheme having the desired features:
\begin{subequations}
\label{CrossRule_2}
\begin{eqnarray}
\epsilon_{-2}^{(n)} & = & \infty \, , 
\qquad \epsilon_{0}^{(n)} \; = \; s_n \, , \qquad n \in \mathbb{N}_0 \, ,
\\
\epsilon_{2k+2}^{(n)} & = & \epsilon_{2k}^{(n+2)} \, + \,
\frac{\displaystyle 
\frac{\Delta \epsilon_{2k}^{(n+1)}}
{\displaystyle \Delta \epsilon_{2k}^{(n)}} - 
\frac{\Delta \epsilon_{2k}^{(n+1)}}
{\displaystyle \epsilon_{2k}^{(n+1)} - \epsilon_{2k-2}^{(n+2)}}}
{\displaystyle 
\frac{1 \rule{0pt}{11pt}}{\displaystyle \Delta \epsilon_{2k}^{(n+1)}} - 
\frac{1}{\displaystyle \Delta \epsilon_{2k}^{(n)}} +
\frac{1}{\displaystyle 
\epsilon_{2k}^{(n+1)} - \epsilon_{2k-2}^{(n+2)}}} \, ,
\qquad k, n \in \mathbb{N}_0 \, .
\end{eqnarray}
\end{subequations}
If we use (\ref{CrossRule_2}) for the computation of
$\epsilon_{2}^{(n)}$, we obtain (\ref{eps_2_2}). 

Next, we use in (\ref{CrossRule_2}) the partial sums (\ref{ParSum_f}) of
the (formal) power series (\ref{PowSer_f}) in the form of
(\ref{ParSum_Rem}). This yields:
\begin{equation}
\epsilon_{2k}^{(n)} \; = \; f (z) \, + \, 
z^{n+2k+1} \, r_{2k}^{(n)} (z) \, , \qquad k, n \in \mathbb{N}_0 \, .
\label{Eps_AccThrOrd}
\end{equation}
The quantities $r_{2k}^{(n)} (z)$ can be computed with the help of the
following recursive scheme which is a direct consequence of the
recursive scheme (\ref{CrossRule_2}) for $\epsilon_{2k}^{(n)}$:
\begin{subequations}
\label{Rec_r}
\begin{eqnarray}
r_{0}^{(n)} (z) & = & 
- \, \sum_{\nu=0}^{\infty} \, \gamma_{n+\nu+1} \, z^{\nu} \; = \;
\frac{f_{n} (z) - f (z)}{z^{n+1}} \, , \qquad n \in \mathbb{N}_0 \, ,
\\
r_{2}^{(n)} (z) & = & r_{0}^{(n+2)} (z) \, + \,
\frac
{\displaystyle \frac{\delta r_{0}^{(n+1)} (z)}{\delta r_{0}^{(n)} (z)}}
{\displaystyle \frac{1 \rule{0pt}{11pt}}{\delta r_{0}^{(n+1)} (z)} - 
\frac{z}{\delta r_{0}^{(n)} (z)}} \, , 
\qquad n \in \mathbb{N}_0 \, ,
\label{Rec_rb}
\\
r_{2k+2}^{(n)} (z) & = & r_{2k}^{(n+2)} (z) \, + \, \frac
{\displaystyle \frac{\delta r_{2k}^{(n+1)} (z)}{\delta r_{2k}^{(n)} (z)} 
- \frac{\delta r_{2k}^{(n+1)} (z)}
{z r_{2k}^{(n+1)} (z) - r_{2k-2}^{(n+2)} (z)}}
{\displaystyle \frac{1 \rule{0pt}{11pt}}{\delta r_{2k}^{(n+1)} (z)} - 
\frac{z}{\delta r_{2k}^{(n)} (z)} + 
\frac{z} {z r_{2k}^{(n+1)} (z) - r_{2k-2}^{(n+2)}(z)}} \, , 
\; k, n \in \mathbb{N}_0 \, . \qquad
\label{Rec_rc}
\end{eqnarray}
\end{subequations}
Here, $\delta r_{2k}^{(n)} (z)$ is defined by (\ref{delta_X}). It should
be noted that (\ref{Rec_rb}) follows from (\ref{Rec_rc}) if we define
$r_{-2}^{(n)} (z) = \infty$.

Similar to the analogous accuracy-through-order relationship
(\ref{Ait_AccThrOrd}) for Aitken's iterated $\Delta^2$ process,
(\ref{Eps_AccThrOrd}) has the right structure to serve as an
accuracy-through-order relationship for Wynn's epsilon algorithm. Thus,
it seems that we have accomplished our aim. However, we are faced with
the same complications as in the case of (\ref{Ait_AccThrOrd}). If
$z^{n+2k+1} r_{2k}^{(n)} (z)$ in (\ref{Eps_AccThrOrd}) is to be of order
$O \bigl( z^{n+2k+1} \bigr)$ as $z \to 0$, then the $z$-independent part
$c_{2k}^{(n)}$ of $r_{2k}^{(n)} (z)$ defined by
\begin{equation}
r_{2k}^{(n)} (z) \; = \; c_{2k}^{(n)} \, + \, O (z) \, , 
\qquad z \to 0 \, ,
\label{Eps_z_indep_NZ_TruncErr_a}
\end{equation} 
has to satisfy
\begin{equation}
c_{2k}^{(n)} \; \ne \; 0 \, , \qquad k, n \in \mathbb{N}_0 \, .
\label{Eps_z_indep_NZ_TruncErr_b}
\end{equation}
If this condition is satisfied, we can be sure that
(\ref{Eps_AccThrOrd}) is indeed the accuracy-through-order relationship
we have been looking for.

As in the case of Aitken's iterated $\Delta^2$ process, it is by no
means obvious whether and how it can be proven that a given power series
gives rise to truncation errors $r_{2k}^{(n)} (z)$ satisfying
(\ref{Eps_z_indep_NZ_TruncErr_a}) and (\ref{Eps_z_indep_NZ_TruncErr_b}). 
Fortunately, it can easily be checked \emph{numerically} whether a given
(formal) power series leads to truncations errors whose $z$-independent
parts are nonzero. If we set $z = 0$ in (\ref{Rec_r}) and use
(\ref{Eps_z_indep_NZ_TruncErr_a}), we obtain the following recursive
scheme:
\begin{subequations}
\label{Rec_c}
\begin{eqnarray}
c_{0}^{(n)} & = & - \, \gamma_{n+1} \, , \qquad n \in \mathbb{N}_0 \, ,
\\
c_{2}^{(n)} & = & c_{0}^{(n+2)} \, - \, 
\frac{\bigl[ c_{0}^{(n+1)} \bigr]^2}{c_{0}^{(n)}} \, ,
\qquad n \in \mathbb{N}_0 \, ,
\label{Rec_c_b}
\\
c_{2k+2}^{(n)} & = & c_{2k}^{(n+2)} \, - \, 
\frac{\bigl[ c_{2k}^{(n+1)} \bigr]^2}{c_{2k}^{(n)}} \, + \, 
\frac{\bigl[ c_{2k}^{(n+1)} \bigr]^2}{c_{2k-2}^{(n+2)}} \, ,
\qquad k \in \mathbb{N} \, , \quad n \in \mathbb{N}_0 \, .
\label{Rec_c_c}
\end{eqnarray}
\end{subequations}
If we define $c_{-2}^{(n)} = \infty$, then (\ref{Rec_c_b}) follows from
(\ref{Rec_c_c}).

Let us now assume that we know for a given (formal) power series that
the $z$-independent parts $c_{2k}^{(n)}$ of the truncation errors
$r_{2k}^{(n)} (z)$ in (\ref{Eps_AccThrOrd}) are nonzero -- either from a
mathematical proof or from a brute force calculation using
(\ref{Rec_c}). Then, (\ref{Eps_AccThrOrd}) is indeed the
accuracy-through-order relationship we have been looking for. This
implies that $\epsilon_{2k}^{(n)}$ can be expressed as follows:
\begin{equation}
\epsilon_{2k}^{(n)} \; = \; f_{n+2k} (z) \, + \, z^{n+2k+1} \, 
\varphi_{2k}^{(n)} (z) \, .
\label{EpsRem_phi}
\end{equation}
If we use this ansatz in (\ref{CrossRule_2}), we obtain the following
recursive scheme:
\begin{subequations}
\label{rec_phi}
\begin{eqnarray}
\varphi_{0}^{(n)} (z) & = & 0 \, , \qquad n \in \mathbb{N}_0 \, ,
\\
\varphi_{2}^{(n)} (z) & = & \frac
{\bigl[ \gamma_{n+2} \bigr]^2}{\gamma_{n+1} - \gamma_{n+2} z} \, ,
\qquad n \in \mathbb{N}_0 \, ,
\label{rec_phi_b}
\\
\varphi_{2k+2}^{(n)} (z) & = & \varphi_{2k}^{(n+2)} (z) \, + \,
\frac{\alpha_{2k+2}^{(n)} (z)}{\beta_{2k+2}^{(n)} (z)} \, ,
\qquad k \in \mathbb{N} \, , \quad n \in \mathbb{N}_0 \, ,
\label{rec_phi_c}
\\
\alpha_{2k+2}^{(n)} (z) & = & 
\frac{\gamma_{n+2k+2} + \delta \varphi_{2k}^{(n+1)} (z)}
{\gamma_{n+2k+1} + \delta \varphi_{2k}^{(n)} (z)} \, - \,
\frac{\gamma_{n+2k+2} + \delta \varphi_{2k}^{(n+1)} (z)}
{\gamma_{n+2k+1} + z \varphi_{2k}^{(n+1)} (z) - 
\varphi_{2k-2}^{(n+2)} (z)} \, ,
\\
\beta_{2k+2}^{(n)} (z) & = & 
\frac{1}{\gamma_{n+2k+2} + \delta \varphi_{2k}^{(n+1)} (z)} \, - \,
\frac{z}{\gamma_{n+2k+1} + \delta \varphi_{2k}^{(n)} (z)} \nonumber \\ 
& & \, + \,
\frac{z}{\gamma_{n+2k+1} + z \varphi_{2k}^{(n+1)} (z) - 
\varphi_{2k-2}^{(n+2)} (z)} \, .
\end{eqnarray}
\end{subequations}
Here, $\delta \varphi_{2k}^{(n)} (z)$ is defined by
(\ref{delta_X}). Moreover, we could also define $\varphi_{-2}^{(n)} (z)
= \infty$. Then, (\ref{rec_phi_b}) would follow from (\ref{rec_phi_c}).

A comparison of (\ref{EpsAsyZero}) and (\ref{EpsRem_phi}) yields
\begin{equation}
\varphi_{2k}^{(n)} (z) \; = \; g_{2k}^{(n)} \, + \, O \bigl( z \bigr) \, ,
\qquad z \to 0 \, .
\label{phi2g}
\end{equation}
Consequently, the $z$-independent part $g_{2k}^{(n)}$ of
$\varphi_{2k}^{(n)} (z)$ is the prediction for the first coefficient
$\gamma_{n+2k+1}$ not used for the computation of $\epsilon_{2k}^{(n)}$.

If we set $z = 0$ in the recursive scheme (\ref{rec_phi}) and use
(\ref{phi2g}), we obtain the following recursive scheme for the
predictions $g_{2k}^{(n)}$:
\begin{subequations}
\label{Rec_g}
\begin{eqnarray}
g_{0}^{(n)} & = & 0 \, , \qquad n \in \mathbb{N}_0 \, ,
\label{Rec_g_a}
\\
g_{2}^{(n)} & = & \frac{\bigl[ \gamma_{n+2} \bigr]^2}{\gamma_{n+1}} \, ,
\qquad n \in \mathbb{N}_0 \, ,
\label{Rec_g_b}
\\
g_{2k+2}^{(n)} & = & g_{2k}^{(n+2)} \, + \, 
\frac
{\bigl[ \gamma_{n+2k+2} - g_{2k}^{(n+1)} \bigr]^2}
{\gamma_{n+2k+1} - g_{2k}^{(n)}} \, - \, 
\frac
{\bigl[ \gamma_{n+2k+2} - g_{2k}^{(n+1)} \bigr]^2}
{\gamma_{n+2k+1} - g_{2k-2}^{(n+2)}} \, , \nonumber \\
& & k \in \mathbb{N} \, , \qquad n \in \mathbb{N}_0 \, .
\label{Rec_g_c}
\end{eqnarray}
\end{subequations}
If we define $g_{-2}^{(n)} = \infty$, then (\ref{Rec_g_b}) follows from
(\ref{Rec_g_a}) and (\ref{Rec_g_c}).

The $z$-independent parts $c_{2k}^{(n)}$ of $r_{2k}^{(n)} (z)$ and
$g_{2k}^{(n)}$ of $\varphi_{2k}^{(n)} (z)$, respectively, are connected.
A comparison of (\ref{Eps_AccThrOrd}), 
(\ref{Eps_z_indep_NZ_TruncErr_a}), (\ref{EpsRem_phi}), and (\ref{phi2g})
yields:
\begin{equation}
g_{2k}^{(n)} \; = \; c_{2k}^{(n)} \, + \, \gamma_{n+2k+1} \, .
\end{equation}

Concerning the choice of the approximation to the limit, we proceed in
the case of the epsilon algorithm just like in the case of Aitken's
iterated $\Delta^2$ process and compute a new approximation to the limit
after the computation of each new partial sum. Thus, if the index $m$ of
the last partial sum $f_{m} (z)$ is even, $m = 2 \mu$, we use as
approximation to the limit $f (z)$ the transformation
\begin{equation}
\bigl\{ f_{0} (z), f_{1} (z), \ldots, f_{2\mu} (z) \bigr\} \mapsto 
\epsilon_{2\mu}^{(0)} \, ,
\end{equation}
and if $m$ is odd, $m = 2\mu+1$, we use the transformation
\begin{equation}
\bigl\{ f_{1} (z), f_{2} (z), \ldots, f_{2\mu+1} (z) \bigr\} \mapsto 
\epsilon_{2\mu}^{(1)} \, .
\end{equation}
These two relationships can be combined into a single equation,
yielding \cite[Eq.\ (4.3-6)]{We89}
\begin{equation}
\bigl\{ f_{m-2\Ent{m/2}} (z), f_{m-2\Ent{m/2}+1} (z), \ldots, 
f_{m} (z) \bigr\} \mapsto \epsilon_{2\Ent{m/2}}^{(m-2\Ent{m/2})} \, ,
\qquad m \in \mathbb{N}_0 \, .
\label{EpsAppLim}
\end{equation}

\setcounter{equation}{0}

\section{The Iteration of Brezinski's Theta Algorithm}
\label{Sec:ThetIt}

Brezinski's theta algorithm is the following recursive scheme
\cite{Bre71}: 
\begin{subequations}
\label{thet_al}
\begin{eqnarray}
\vartheta_{-1}^{(n)} & = & 0 \, , \qquad
\vartheta_0^{(n)} \; = \; s_n \, , \qquad n \in \mathbb{N}_0 \, , 
\\
\vartheta_{2 k + 1}^{(n)} & = & \vartheta_{2 k-1}^{(n+1)}
\, + \, 1 / \bigl[\Delta \vartheta_{2 k}^{(n)}\bigr] \, , 
\qquad k, n \in \mathbb{N}_0 \, ,
\label{thet_al_b}
\\
\vartheta_{2 k+2}^{(n)} & = & \vartheta_{2 k}^{(n+1)} \, +
\, \frac
{\bigl[ \Delta \vartheta_{2 k}^{(n+1)} \bigr] \,
\bigl[\Delta \vartheta_{2 k + 1}^{(n+1)} \bigr]}
{\Delta^2 \vartheta_{2 k+1}^{(n)}} \, , 
\qquad k, n \in \mathbb{N}_0 \, .
\end{eqnarray}
\end{subequations}
As in the case of Wynn's epsilon algorithm (\ref{eps_al}), only the
elements $\vartheta_{2k}^{(n)}$ with even subscripts provide
approximations to the (generalized) limit of the sequence to be
transformed. The elements $\vartheta_{2k+1}^{(n)}$ with odd subscripts
are only auxiliary quantities which diverge if the whole process
converges. 

The theta algorithm was derived from Wynn's epsilon algorithm
(\ref{eps_al}) with the intention of overcoming the inability of the
epsilon algorithm to accelerate logarithmic convergence. In that
respect, the theta algorithm was a great success. Extensive numerical
studies of Smith and Ford \cite{SmiFor79,SmiFor82} showed that the theta
algorithm is not only very powerful, but also much more versatile than
the epsilon algorithm. Like the epsilon algorithm, it is an efficient
accelerator for linear convergence and it is also able to sum many
divergent series. However, it is also able to accelerate the convergence
of many logarithmically convergent sequences and series. 

As for example discussed in \cite{We91}, new sequence transformations
can be constructed by iterating explicit expressions for sequence
transformations with low transformation orders. The best known example
of such an iterated sequence transformation is probably Aitken's
iterated $\Delta^2$ process (\ref{ItAit_1}) which is obtained by
iterating Aitken's $\Delta^2$ formula (\ref{AitFor_1}).

The same approach is also possible in the case of the theta algorithm. A
suitable closed-form expression, which may be iterated, is \cite[Eq.\
(10.3-1)]{We89}
\begin{equation}
\vartheta_2^{(n)} \; = \; s_{n+1} \, - \, \frac
{\bigl[\Delta s_n\bigr] \bigl[\Delta s_{n+1}\bigr] 
\bigl[\Delta^2 s_{n+1}\bigr]}
{\bigl[\Delta s_{n+2}\bigr] \bigl[\Delta^2 s_n\bigr] -
\bigl[\Delta s_n\bigr] \bigl[\Delta^2 s_{n+1}\bigr]} \, , 
\qquad n \in \mathbb{N}_0 \, .
\label{theta_2_1}
\end{equation}
The iteration of this expression yields the following nonlinear
recursive scheme \cite[Eq.\ (10.3-6)]{We89}:
\begin{subequations}
\label{thetit_1}
\begin{eqnarray}
\mathcal{J}_0^{(n)} & = & s_n \, , \qquad n \in \mathbb{N}_0 \, ,
\\
\mathcal{J}_{k+1}^{(n)} & = &
\mathcal{J}_k^{(n+1)} \, - \, \frac
{\bigl[ \Delta \mathcal{J}_k^{(n)} \bigr] 
\bigl[ \Delta \mathcal{J}_k^{(n+1)} \bigr]
\bigl[ \Delta^2 \mathcal{J}_k^{(n+1)} \bigr]}
{\bigl[ \Delta \mathcal{J}_k^{(n+2)} \bigr] 
\bigl[ \Delta^2 \mathcal{J}_k^{(n)} \bigr]  - 
\bigl[ \Delta \mathcal{J}_k^{(n)} \bigr] 
\bigl[\Delta^2 \mathcal{J}_k^{(n+1)} \bigr]}
\, , \qquad k,n \in \mathbb{N}_0 \, .
 \end{eqnarray}
\end{subequations}
In convergence acceleration and summation processes, the iterated
transformation $\mathcal{J}_k^{(n)}$ has similar properties as the theta
algorithm from which it was derived: They are both very powerful as well
as very versatile. $\mathcal{J}_k^{(n)}$ is not only an effective
accelerator for linear convergence as well as able to sum divergent
series, but it is also able to accelerate the convergence of many
logarithmically convergent sequences and series
\cite{BhoBhaRoy89,Sab87,Sab91,Sab92,Sab95,We89,We91,We94b}.

In spite of all these similarities, the iterated transformation
$\mathcal{J}_{k}^{(n)}$ has one undeniable advantage over the theta
algorithm, which ultimately explains why in this article only
$\mathcal{J}_{k}^{(n)}$ is studied, but not the theta algorithm: The
recursive scheme (\ref{thetit_1}) for $\mathcal{J}_{k}^{(n)}$ is
slightly less complicated than the recursive scheme (\ref{thet_al}) for
the theta algorithm. On p.\ 282 of \cite{We89} it was emphasized that a
replacement of (\ref{thet_al_b}) by the simpler recursion
\begin{equation}
\vartheta_{2 k + 1}^{(n)} \; = \;
1 / \bigl[\Delta \vartheta_{2 k}^{(n)}\bigr] \, , 
\qquad k, n \in \mathbb{N}_0 \, ,  
\end{equation}
would lead to a modified theta algorithm which satisfies $\vartheta_{2
k}^{(n)} = {\cal J}_k^{(n)}$.

It is a direct consequence of the recursive scheme (\ref{thetit_1}) that
$3k+1$ sequence elements $s_{n}$, $s_{n+1}$, \ldots, $s_{n+3k}$ are
needed for the computation of $\mathcal{J}_{k}^{(n)}$. Thus, we now
choose as input data the partial sums (\ref{ParSum_f}) of the (formal)
power series (\ref{PowSer_f}) according to $s_n = f_n (z)$, and
conjecture that all coefficients $\gamma_{0}$, $\gamma_{1}$, \ldots,
$\gamma_{n+3k}$, which were used for the construction of
$\mathcal{J}_{k}^{(n)}$, are exactly reproduced by a Taylor
expansion. This means that we have to look for an accuracy-through-order
relationship of the following kind:
\begin{equation}
f (z) \, - \, \mathcal{J}_{k}^{(n)} \; = \; 
O \bigl( z^{n+3k+1} \bigr) \, , \qquad z \to 0 \, .
\label{tit_OrdEst}
\end{equation}
Such an accuracy-through-order relationship would imply that
$\mathcal{J}_{k}^{(n)}$ can be expressed as follows:
\begin{equation}
\mathcal{J}_{k}^{(n)} \; = \; f_{n+3k} (z) \, + \, 
\mathcal{G}_{k}^{(n)} \, z^{n+3k+1} \, + \, 
O \bigl( z^{n+3k+2} \bigr) \, , \qquad z \to 0 \, .
\label{tit_AsyZero}
\end{equation}
The constant $\mathcal{G}_{k}^{(n)}$ is the prediction made for the
coefficient $\gamma_{n+3k+1}$, which is the first coefficient of the
power series (\ref{PowSer_f}) not used for the computation of
$\mathcal{J}_{k}^{(n)}$.

Unfortunately, the recursive scheme (\ref{thetit_1}) is not suited for
our purposes. This can be shown by computing $\mathcal{J}_{1}^{(n)}$
from the partial sums $f_{n} (z)$, $f_{n+1} (z)$, $f_{n+2} (z)$, and
$f_{n+3} (z)$:
\begin{equation}
\mathcal{J}_{1}^{(n)} \; = \; f_{n+1} (z) \, - \, \frac
{\gamma_{n+1} \gamma_{n+2} [\gamma_{n+3} z - \gamma_{n+2}] z^{n+2}}
{\gamma_{n+3} z [\gamma_{n+2} z - \gamma_{n+1}] - \gamma_{n+1} 
[\gamma_{n+3} z - \gamma_{n+2}]} \, .
\end{equation}
Superficially, it looks as if the accuracy-through-order relationship
(\ref{tit_OrdEst}) is not satisfied by $\mathcal{J}_{1}^{(n)}$. However,
the rational expression on the right-hand side contains the missing
terms $\gamma_{n+2} z^{n+2}$ and $\gamma_{n+3} z^{n+3}$, as shown by the
Taylor expansion
\begin{eqnarray}
\lefteqn{- \, \frac
{\gamma_{n+1} \gamma_{n+2} [\gamma_{n+3} z - \gamma_{n+2}] z^{n+2}}
{\gamma_{n+3} z [\gamma_{n+2} z - \gamma_{n+1}] - \gamma_{n+1} 
[\gamma_{n+3} z - \gamma_{n+2}]}} \nonumber \\
& \qquad = & \gamma_{n+2} z^{n+2} \, + \, \gamma_{n+3} z^{n+3} \, - \, 
\frac{\gamma_{n+3} 
\Bigl\{ \bigl[ \gamma_{n+2} \bigr]^2 - 
2 \gamma_{n+1} \gamma_{n+3} \Bigr\} z^{n+4}}
{\gamma_{n+1} \gamma_{n+2}} \, + \, 
O \bigl( z^{n+5} \bigr) \, . \quad
\label{rem_theta_2}
\end{eqnarray}
Thus, an expression, which is in agreement with (\ref{tit_AsyZero}), can
be obtained easily in the case of the simplest transform
$\mathcal{J}_{1}^{(n)}$. Moreover, the Taylor expansion
(\ref{rem_theta_2}) shows that $\mathcal{J}_{1}^{(n)}$ makes the
prediction
\begin{equation}
\mathcal{G}_{1}^{(n)} \; = \; - \, \frac{\gamma_{n+3} 
\Bigl\{ \bigl[ \gamma_{n+2} \bigr]^2 - 
2 \gamma_{n+1} \gamma_{n+3} \Bigr\}}
{\gamma_{n+1} \gamma_{n+2}}
\end{equation}
for the first series coefficient $\gamma_{n+4}$ not used for the
computation of $\mathcal{J}_{1}^{(n)}$. Of course, by including
additional terms in the Taylor expansion (\ref{rem_theta_2}) further
predictions on series coefficients with higher indices can be made.

However, in the case of more complicated transforms
$\mathcal{J}_{k}^{(n)}$ with $k > 1$ it by no means obvious whether and
how an expression, which is in agreement with (\ref{tit_AsyZero}), can
be constructed. Consequently, it is certainly a good idea to replace the
recursive scheme (\ref{thetit_1}) by an alternative recursive scheme,
which directly leads to appropriate expressions for
$\mathcal{J}_{k}^{(n)}$ with $k > 1$.

Many different expressions for $\vartheta_{2}^{(n)}$ in terms of
$s_n$, $s_{n+1}$, $s_{n+2}$, and $s_{n+3}$ are known 
\cite[Section 10.4]{We89}. The for our purposes appropriate expression 
is
\begin{equation}
\vartheta_{2}^{(n)} \; = \; s_{n+3} \, - \, 
\frac{\bigl[ \Delta s_{n+2} \bigr] \Bigl\{\bigl[ \Delta s_{n+2} \bigr] 
\bigl[\Delta^2 s_{n} \bigr] + \bigl[ \Delta s_{n+1} \bigr]^2 - 
\bigl[ \Delta s_{n+2} \bigr] \bigl[ \Delta s_{n} \bigr] \Bigr\}}
{\bigl[ \Delta s_{n+2} \bigr] \bigl[ \Delta^2 s_{n} \bigr] - 
\bigl[ \Delta s_{n} \bigr] \bigl[ \Delta^2 s_{n+1} \bigr]} \, .
\label{theta_2_2}
\end{equation}
Just like (\ref{theta_2_1}), this expression can be iterated and yields
\begin{subequations}
\label{thetit_2}
\begin{eqnarray}
\mathcal{J}_0^{(n)} & = & s_n \, , \qquad n \in \mathbb{N}_0 \, ,
\\
\mathcal{J}_{k+1}^{(n)} & = & \mathcal{J}_{k}^{(n+3)} \, - \,
\frac{A_{k+1}^{(n)}}{B_{k+1}^{(n)}} \, ,
\qquad k, n \in \mathbb{N}_0 \, ,
\\
A_{k+1}^{(n)} & = & \bigl[ \Delta \mathcal{J}_{k}^{(n+2)} \bigr]
\Bigl\{ \bigl[ \Delta \mathcal{J}_{k}^{(n+2)} \bigr] 
\bigl[ \Delta^2 \mathcal{J}_{k}^{(n)} \bigr] \, + \,
\bigl[ \Delta \mathcal{J}_{k}^{(n+1)} \bigr]^2 \nonumber \\
& & \quad \, - \,
\bigl[ \Delta \mathcal{J}_{k}^{(n)} \bigr] 
\bigl[ \Delta \mathcal{J}_{k}^{(n+2)} \bigr] \Bigr\} \, ,
\\
B_{k+1}^{(n)} & = & \bigl[ \Delta \mathcal{J}_{k}^{(n+2)} \bigr]
\bigl[ \Delta^2 \mathcal{J}_{k}^{(n)} \bigr] \, - \,
\bigl[ \Delta \mathcal{J}_{k}^{(n)} \bigr]
\bigl[ \Delta^2 \mathcal{J}_{k}^{(n+1)} \bigr] \, .
\end{eqnarray}
\end{subequations}
If we now use either (\ref{theta_2_2}) or (\ref{thetit_2}) to compute
$\mathcal{J}_{1}^{(n)}$ from the partial sums $f_{n} (z)$, $f_{n+1}
(z)$, $f_{n+2} (z)$, and $f_{n+3} (z)$, we obtain the following
expression which obviously possesses the desired features:
\begin{equation}
\mathcal{J}_{1}^{(n)} \; = \; f_{n+3} (z) \, - \, \frac
{\gamma_{n+3} \Bigl\{ \gamma_{n+3} 
\bigl[ \gamma_{n+2} z - \gamma_{n+1} \bigr] + 
\bigl[ \gamma_{n+2} \bigr]^2 - \gamma_{n+1} \gamma_{n+3} \Bigr\}
z^{n+4}}
{\gamma_{n+3} z \bigl[ \gamma_{n+2} z - \gamma_{n+1} \bigr] -
\gamma_{n+1} \bigl[ \gamma_{n+3} z - \gamma_{n+2} \bigr]} \, .
\end{equation}

Next, we use in (\ref{thetit_2}) the partial sums (\ref{ParSum_f}) of
the (formal) power series (\ref{PowSer_f}) in the form of
(\ref{ParSum_Rem}). This yields:
\begin{equation}
\mathcal{J}_{k}^{(n)} \; = \; f (z) \, + \, 
z^{n+3k+1} \, \mathcal{R}_{k}^{(n)} (z) \, , 
\qquad k, n \in \mathbb{N}_0 \, .
\label{tit_AccThrOrd}
\end{equation}
The quantities $\mathcal{R}_{k}^{(n)} (z)$ can be computed with the help
of the following recursive scheme which is a direct consequence of the
recursive scheme (\ref{thetit_2}) for $\mathcal{J}_{k}^{(n)}$:
\begin{subequations}
\label{Rec_calR}
\begin{eqnarray}
\mathcal{R}_{0}^{(n)} (z) & = & 
- \, \sum_{\nu=0}^{\infty} \, \gamma_{n+\nu+1} \, z^{\nu} \; = \;
\frac{f_{n} (z) - f (z)}{z^{n+1}} \, , \qquad n \in \mathbb{N}_0 \, ,
\\
\mathcal{R}_{k+1}^{(n)} (z) & = & \mathcal{R}_{k}^{(n+3)} (z) \, - \,
\frac{\mathcal{N}_{k+1}^{(n)} (z)}{\mathcal{D}_{k+1}^{(n)} (z)} \, , 
\qquad k, n \in \mathbb{N}_0 \, ,
\\
\mathcal{N}_{k+1}^{(n)} (z) & = & 
\bigl[ \delta \mathcal{R}_{k}^{(n+2)} (z) \bigr]
\Bigl\{ \bigl[ \delta \mathcal{R}_{k}^{(n+2)} (z) \bigr]
\bigl[ \delta^2 \mathcal{R}_{k}^{(n)} (z) \bigr] \, + \,
\bigl[ \delta \mathcal{R}_{k}^{(n+1)} (z) \bigr]^2 \nonumber \\ 
& & \quad \, - \,
\bigl[ \delta \mathcal{R}_{k}^{(n)} (z) \bigr]
\bigl[ \delta \mathcal{R}_{k}^{(n+2)} (z) \bigr] \Bigr\} \, ,
\\
\mathcal{D}_{k+1}^{(n)} (z) & = &
z \bigl[ \delta \mathcal{R}_{k}^{(n+2)} (z) \bigr] 
\bigl[ \delta^2 \mathcal{R}_{k}^{(n)} (z) \bigr] \, - \,
\bigl[ \delta \mathcal{R}_{k}^{(n)} (z) \bigr]
\bigl[ \delta^2 \mathcal{R}_{k}^{(n+1)} (z) \bigr] \, .
\end{eqnarray}
\end{subequations}
Here, $\delta \mathcal{R}_{k}^{(n+2)} (z)$ and $\delta^2
\mathcal{R}_{k}^{(n+2)} (z)$ are defined by (\ref{delta_X}).

Similar to the analogous accuracy-through-order relationships
(\ref{Ait_AccThrOrd}) and (\ref{Eps_AccThrOrd}) for Aitken's iterated
$\Delta^2$ process and the epsilon algorithm, respectively,
(\ref{tit_AccThrOrd}) has the right structure to serve as an
accuracy-through-order relationship for the iterated theta algorithm.
Thus, it seems that we have accomplished our aim. However, we are faced
with the same complications as in the case of (\ref{Ait_AccThrOrd}) and
(\ref{Eps_AccThrOrd}). If $z^{n+3k+1} \mathcal{R}_{2k}^{(n)} (z)$ in
(\ref{tit_AccThrOrd}) is to be of order $O \bigl( z^{n+3k+1} \bigr)$ as
$z \to 0$, then the $z$-independent part $\mathcal{C}_{k}^{(n)}$ of
$\mathcal{R}_{k}^{(n)} (z)$ defined by
\begin{equation}
\mathcal{R}_{k}^{(n)} (z) \; = \; \mathcal{C}_{k}^{(n)} \, + \, O (z) \, , 
\qquad z \to 0 \, ,
\label{tit_z_indep_NZ_TruncErr_a}
\end{equation} 
has to satisfy
\begin{equation}
\mathcal{C}_{k}^{(n)} \; \ne \; 0 \, , \qquad k, n \in \mathbb{N}_0 \, .
\label{tit_z_indep_NZ_TruncErr_b}
\end{equation}
If this condition is satisfied, then it is guaranteed that
(\ref{tit_AccThrOrd}) is indeed the accuracy-through-order relationship
we have been looking for.

As in the case of Aitken's iterated $\Delta^2$ process or the epsilon
algorithm, it is by no means obvious whether and how it can be proven
that a given power series gives rise to truncation errors
$\mathcal{R}_{k}^{(n)} (z)$ satisfying (\ref{tit_z_indep_NZ_TruncErr_a})
and (\ref{tit_z_indep_NZ_TruncErr_b}). Fortunately, it can easily be
checked \emph{numerically} whether a given (formal) power series leads
to truncations errors whose $z$-independent parts are nonzero. If we set
$z = 0$ in (\ref{Rec_calR}) and use (\ref{tit_z_indep_NZ_TruncErr_a}),
we obtain the following recursive scheme:
\begin{subequations}
\label{Rec_calG}
\begin{eqnarray}
\mathcal{C}_{0}^{(n)} & = & - \, \gamma_{n+1} \, , 
\qquad n \in \mathbb{N}_0 \, ,
\\
\mathcal{C}_{k+1}^{(n)} & = & \mathcal{C}_{k}^{(n+3)} \, - \, \frac
{\mathcal{C}_{k}^{(n+2)} 
\Bigl\{2 \mathcal{C}_{k}^{(n)} \mathcal{C}_{k}^{(n+2)} - 
\bigl[ \mathcal{C}_{k}^{(n+1)} \bigr]^2 \Bigr\}}
{\mathcal{C}_{k}^{(n)} \mathcal{C}_{k}^{(n+1)}} \, , 
\qquad k, n \in \mathbb{N}_0 \, .
\end{eqnarray}
\end{subequations}

Let us now assume that we know for a given (formal) power series that
the $z$-independent parts $\mathcal{G}_{k}^{(n)}$ of the truncation
errors $\mathcal{R}_{k}^{(n)} (z)$ in (\ref{tit_AccThrOrd}) are nonzero
-- either from a mathematical proof or from a brute force calculation
using (\ref{Rec_calG}). Then, (\ref{tit_AccThrOrd}) is indeed the
accuracy-through-order relationship we have been looking for. This
implies that $\mathcal{J}_{k}^{(n)}$ can be expressed as follows:
\begin{equation}
\mathcal{J}_{k}^{(n)} \; = \; f_{n+3k} (z) \, + \, 
z^{n+3k+1} \, \Psi_{k}^{(n)} (z) \, , 
\qquad k, n \in \mathbb{N}_0 \, .
\label{titRemPsi}
\end{equation}
If we use this ansatz in (\ref{thetit_2}), we obtain the following
recursive scheme:
\begin{subequations}
\label{Rec_Psi}
\begin{eqnarray}
\Psi_{0}^{(n)} (z) & = & 0 \, , \qquad n \in \mathbb{N}_0 \, ,
\\
\Psi_{1}^{(n)} (z) & = & - \, 
\frac{\gamma_{n+3} 
\Bigl\{ \gamma_{n+3} \bigl[ \gamma_{n+2} z - \gamma_{n+1} \bigr] + 
\bigl[ \gamma_{n+2} \bigr]^2 - \gamma_{n+1} \gamma_{n+3} \Bigr\}}
{\gamma_{n+3} \bigl[ \gamma_{n+2} z - \gamma_{n+1} \bigr] - 
\gamma_{n+1} \bigl[ \gamma_{n+3} z - \gamma_{n+2} \bigr]}
\, , \qquad n \in \mathbb{N}_0 \, , \qquad
\\
\Psi_{k+1}^{(n)} (z) & = & \Psi_{k}^{(n+3)} (z) \, - \,
\frac{N_{k+1}^{(n)} (z)}{D_{k+1}^{(n)} (z)}
\, , \qquad k, n \in \mathbb{N}_0 \, ,
\\
N_{k+1}^{(n)} (z) & = & 
\bigl[ \gamma_{n+3k+3} + \delta \Psi_{k}^{(n+2)} (z) \bigr] \nonumber \\
& & \quad \times \, 
\Bigl\{ \bigl[ \gamma_{n+3k+3} + \delta \Psi_{k}^{(n+2)} (z) \bigr]
\bigl[ \gamma_{n+3k+2} z - \gamma_{n+3k+1} + 
\delta^2 \Psi_{k}^{(n)} (z) \bigr] \nonumber \\ 
& & \quad \qquad + \,
\bigl[ \gamma_{n+3k+2} + \delta \Psi_{k}^{(n+1)} (z) \bigr]^2 \nonumber
\\
& & \qquad \qquad \quad - \,
\bigl[ \gamma_{n+3k+1} + \delta \Psi_{k}^{(n)} (z) \bigr]
\bigl[ \gamma_{n+3k+3} + \delta \Psi_{k}^{(n+2)} (z) \bigr] \Bigr\} \, ,
\\
D_{k+1}^{(n)} (z) & = & 
\bigl[ \gamma_{n+3k+3} + \delta \Psi_{k}^{(n+2)} (z) \bigr]
\bigl[ \gamma_{n+3k+2} z - \gamma_{n+3k+1} + 
\delta^2 \Psi_{k}^{(n)} (z) \bigr] \nonumber \\
& & \qquad - \,
\bigl[ \gamma_{n+3k+1} + \delta \Psi_{k}^{(n)} (z) \bigr]
\bigl[ \gamma_{n+3k+3} z - \gamma_{n+3k+2} + 
\delta^2 \Psi_{k}^{(n+1)} (z) \bigr] \, .
\end{eqnarray}
\end{subequations}
Here, $\delta \Psi_{k}^{(n+2)} (z)$ and $\delta^2
\Psi_{k}^{(n+2)} (z)$ are defined by (\ref{delta_X}).

A comparison of (\ref{tit_AsyZero}) and (\ref{titRemPsi}) yields
\begin{equation}
\Psi_{k}^{(n)} (z) \; = \; \mathcal{G}_{k}^{(n)} \, + \, 
O \bigl( z \bigr) \, , \qquad z \to 0 \, .
\label{Psi2calG}
\end{equation}
Consequently, the $z$-independent part $\mathcal{G}_{k}^{(n)}$ of
$\Psi_{k}^{(n)} (z)$ is the prediction for the first coefficient
$\gamma_{n+3k+1}$ not used for the computation of
$\mathcal{J}_{k}^{(n)}$.

If we set $z = 0$ in the recursive scheme (\ref{Rec_Psi}) and use
(\ref{Psi2calG}), we obtain the following recursive scheme for the
predictions $\mathcal{G}_{k}^{(n)}$:
\begin{subequations}
\label{Rec_mclG}
\begin{eqnarray}
\mathcal{G}_{0}^{(n)} & = & 0 \, , \qquad n \in \mathbb{N}_0 \, ,
\\
\mathcal{G}_{1}^{(n)} & = & - \, \frac{\gamma_{n+3} 
\Bigl\{ \bigl[ \gamma_{n+2} \bigr]^2 - 
2 \gamma_{n+1} \gamma_{n+3} \Bigr\}}
{\gamma_{n+1} \gamma_{n+2}}
\, , \qquad n \in \mathbb{N}_0 \, ,
\\
\mathcal{G}_{k+1}^{(n)} & = & \mathcal{G}_{k}^{(n+3)} \, - \, 
\frac{F_{k+1}^{(n)}}{H_{k+1}^{(n)}}
\, , \qquad k, n \in \mathbb{N}_0 \, ,
\\
F_{k+1}^{(n)} & = & 
\bigl[ \gamma_{n+3k+3} - \mathcal{G}_{k}^{(n+2)} \bigr]
\Bigl\{ \bigl[ \gamma_{n+3k+2} - \mathcal{G}_{k}^{(n+1)} \bigr]^2
\nonumber \\
& & \qquad - \, 
2 \, \bigl[ \gamma_{n+3k+1} - \mathcal{G}_{k}^{(n)} \bigr]
\bigl[ \gamma_{n+3k+3} - \mathcal{G}_{k}^{(n+2)} \bigr] \Bigr\}
\, ,
\\
H_{k+1}^{(n)} & = &
\bigl[ \gamma_{n+3k+1} - \mathcal{G}_{k}^{(n)} \bigr]
\bigl[ \gamma_{n+3k+2} - \mathcal{G}_{k}^{(n+1)} \bigr] \, .
\end{eqnarray}
\end{subequations}

The $z$-independent parts $\mathcal{C}_{k}^{(n)}$ of
$\mathcal{R}_{k}^{(n)} (z)$ and $\mathcal{G}_{k}^{(n)}$ of
$\Psi_{k}^{(n)} (z)$, respectively, are connected. A comparison of
(\ref{tit_AccThrOrd}), (\ref{tit_z_indep_NZ_TruncErr_a}),
(\ref{titRemPsi}), and (\ref{Psi2calG}) yields:
\begin{equation}
\mathcal{G}_{k}^{(n)} \; = \; \mathcal{C}_{k}^{(n)} \, + \, 
\gamma_{n+3k+1} \, .
\end{equation}

As in the case of Aitken's iterated $\Delta^2$ process or Wynn's epsilon 
algorithm, a new approximation to the limit will be computed after the
computation of each new partial sum. Thus, if the index $m$ of
the last partial sum $f_{m} (z)$ is a multiple of 3, $m = 3\mu$, we use as
approximation to the limit $f (z)$ the transformation
\begin{equation}
\bigl\{ f_{0} (z), f_{1} (z), \ldots, f_{3\mu} (z) \bigr\} \mapsto 
\mathcal{}J_{\mu}^{(0)} \, ,
\end{equation}
if we have $m = 3\mu+1$, we use the transformation
\begin{equation}
\bigl\{ f_{1} (z), f_{2} (z), \ldots, f_{3\mu+1} (z) \bigr\} \mapsto 
\mathcal{J}_{\mu}^{(1)} \, ,
\end{equation}
and if we have $m = 3\mu+2$, we use the transformation
\begin{equation}
\bigl\{ f_{2} (z), f_{3} (z), \ldots, f_{3\mu+2} (z) \bigr\} \mapsto 
\mathcal{J}_{\mu}^{(2)} \, ,
\end{equation}
These three relationships can be combined into a single equation,
yielding \cite[Eq.\ (10.4-7)]{We89}
\begin{equation}
\bigl\{ f_{m-3\Ent{m/3}} (z), f_{m-3\Ent{m/3}+1} (z), \ldots, 
f_{m} (z) \bigr\} \mapsto \mathcal{J}_{\Ent{m/3}}^{(m-3\Ent{m/3})} \, ,
\qquad m \in \mathbb{N}_0 \, .
\label{TitAppLim}
\end{equation}

\setcounter{equation}{0}

\section{Applications}
\label{Sec:Applics}

In this article, two principally different kinds of results were
derived. The first group of results -- the accuracy-through-order
relationships (\ref{Ait_AccThrOrd}), (\ref{Eps_AccThrOrd}), and
(\ref{tit_AccThrOrd}) and the corresponding recursive schemes
(\ref{RecR}), (\ref{Eps_AccThrOrd}), and (\ref{Rec_calR}) -- defines the
transformation error terms $z^{n+2k+1} R_{k}^{(n)} (z)$, $z^{n+2k+1}
r_{2k}^{(n)} (z)$, and $z^{n+3k+1} \mathcal{R}_{k}^{(n)} (z)$. These
quantities describe how the rational approximants
$\mathcal{A}_{k}^{(n)}$, $\epsilon_{2k}^{(n)}$, and
$\mathcal{J}_{k}^{(n)}$ differ from the function $f (z)$ which is to be
approximated. Obviously, the transformation error terms must vanish if
the transformation process converges.

The second group of results -- (\ref{AitRemPhi}), (\ref{EpsRem_phi}),
and (\ref{titRemPsi}) and the corresponding recursive schemes
(\ref{RecPhi}), (\ref{rec_phi}), and (\ref{Rec_Psi}) -- defines the
terms $z^{n+2k+1} \Phi_{k}^{(n)} (z)$, $z^{n+2k+1} 
\varphi_{2k}^{(n)} (z)$, and $z^{n+3k+1} \Psi_{k}^{(n)} (z)$. These 
quantities describe how the rational approximants
$\mathcal{A}_{k}^{(n)}$, $\epsilon_{2k}^{(n)}$, and
$\mathcal{J}_{k}^{(n)}$ differ from the partial sums $f_{n+2k} (z)$ and 
$f_{n+3k} (z)$, respectively, from which they were constructed. Hence,
the first group of results essentially describes what is still missing
in the transformation process, whereas the second group describes what
was gained by constructing rational expressions from the partial sums.

The recursive schemes (\ref{RecR}), (\ref{Eps_AccThrOrd}), and
(\ref{Rec_calR}) of the first group use as input data the remainder
terms
\begin{equation}
\frac{f_{n} (z) - f (z)}{z^{n+1}} \; = \; - \, 
\sum_{\nu=0}^{\infty} \, \gamma_{n+\nu+1} \, z^{\nu} \, .
\label{f_n_Rem}
\end{equation}
In most practically relevant convergence acceleration and summation
problems, only a finite number of series coefficients $\gamma_{\nu}$ are
known. Consequently, the remainder terms (\ref{f_n_Rem}) are usually not
known explicitly, which means that the immediate practical usefulness of
the first group of results is quite limited.  Nevertheless, these
results are of interest because they can be used to study the
convergence of the sequence transformations of this article for model
problems.

As an example, let us consider the following series expansion for the
logarithm, 
\begin{equation}
\frac{\ln(1+z)}{z} \; = \; {}_2 F_1 (1, 1; 2; - z) \; = \; 
\sum_{m=0}^{\infty} \, \frac{(-z)^m}{m+1} \, ,
\label{LogSer}
\end{equation}
which converges for all $z \in \mathbb{C}$ with $\vert z \vert < 1$. The 
logarithm possesses the integral representation
\begin{equation}
\frac{\ln(1+z)}{z} \; = \; \int_{0}^{1} \, \frac{\mathrm{d} t}{1+zt} \, ,
\end{equation}
which shows that $\ln(1+z)/z$ is a Stieltjes function and that the
hypergeometric series on the right-hand side of (\ref{LogSer}) is the
corresponding Stieltjes series (a detailed treatment of Stieltjes
functions and Stieltjes series can for example be found in Section 5 of
\cite{BaGM96}). Consequently, $\ln(1+z)/z$ possesses the following
representation as a partial sum plus an explicit remainder which is
given by a Stieltjes integral (compare for example Eq.\ (13.1-5) of
\cite{We89}): 
\begin{equation}
\frac{\ln(1+z)}{z} \; = \; \sum_{\nu=0}^{n} \, \frac{(-z)^m}{m+1} 
\, + \, (-z)^{n+1} \, 
\int_{0}^{1} \, \frac{t^{n+1} \mathrm{d} t}{1+zt} \, , 
\qquad n \in \mathbb{N}_0 \, .
\end{equation}
For $\vert z \vert < 1$, the numerator of the remainder integral on the
right-hand side can be expanded. Interchanging summation and integration
then yields:
\begin{equation}
(-1)^{n+1} \, \int_{0}^{1} \, \frac{t^{n+1} \mathrm{d} t}{1+zt} \; = \; 
\sum_{m=0}^{\infty} \, \frac{(-1)^{n+m+1} z^m}{n+m+2} \, .
\label{ln_RemExpr}
\end{equation}
Next, we use for $0 \le n \le 6$ the negative of these remainder
integrals as input data in the recursive schemes (\ref{RecR}),
(\ref{Eps_AccThrOrd}), and (\ref{Rec_calR}), and do a Taylor expansion
of the resulting expressions. Thus, we obtain according to
(\ref{Ait_AccThrOrd}), (\ref{Eps_AccThrOrd}), and (\ref{tit_AccThrOrd}):
\begin{subequations}
\label{f_z_Pred}
\begin{eqnarray}
\mathcal{A}_{3}^{(0)} & = & \frac{\ln(1+z)}{z} \, + \, 
{\frac {421 z^7}{16537500}} \, - \, 
{\frac {796321 z^8}{8682187500}} \, + \, 
{\frac {810757427 z^9}{4051687500000}} \, + \, 
O \bigl( z^{10} \bigr) \, ,
\\
\epsilon_{6}^{(0)} & = & \frac{\ln(1+z)}{z} \, + \, 
{\frac {z^7}{9800}} \, - \, {\frac {31 z^8}{77175}} 
\, + \, {\frac {113 z^9}{120050}}
\, + \, O \bigl( z^{10} \bigr) \, ,
\\
\mathcal{J}_{2}^{(0)} & = & \frac{\ln(1+z)}{z} \, + \, 
{\frac {z^7}{37800}} \, - \, {\frac {19 z^8}{198450}}
\, + \, {\frac {z^9}{4725}} 
\, + \, O \bigl( z^{10} \bigr) \, .
\end{eqnarray}
\end{subequations}
All calculations were done symbolically, using the exact rational
arithmetics of Maple. Consequently, the results in (\ref{f_z_Pred}) are
exact and free of rounding errors.

The leading coefficients of the Taylor expansions of the transformation
error terms for $\mathcal{A}_{3}^{(0)}$ and $\mathcal{J}_{2}^{(0)}$ are
evidently smaller than the corresponding coefficients for
$\epsilon_{6}^{(0)}$. This observation provides considerable evidence
that Aitken's iterated $\Delta^2$ process and Brezinski's iterated theta 
algorithm are in the case of the series (\ref{LogSer}) for $\ln(1+z)/z$
more effective than Wynn's epsilon algorithm which according to
(\ref{Eps_Pade}) produces Pad\'{e} approximants. 

This conclusion is also confirmed by the following numerical example in
Table I, in which the convergence of the series (\ref{LogSer}) for
$\ln(1+z)/z$ is accelerated for $z = 0.95$. The numerical values of the
remainder terms (\ref{ln_RemExpr}) were used as input data in the
recursive schemes (\ref{RecR}), (\ref{Eps_AccThrOrd}), and
(\ref{Rec_calR}) to compute numerically the transformation error terms
in (\ref{Ait_AccThrOrd}), (\ref{Eps_AccThrOrd}), and
(\ref{tit_AccThrOrd}). The transformation error terms, which are listed
in columns 3 - 5, were chosen in agreement with (\ref{AitAppLim}),
(\ref{EpsAppLim}), and (\ref{TitAppLim}), respectively.

\begin{table}[htb]
\begin{tabular}{lllll}
\multicolumn{5}{l}{\bf Table I: Convergence of the Transformation Error 
Terms} \\
\multicolumn{5}{l}{Transformation of 
$\ln(1+z)/z = \sum_{m=0}^{\infty} (-z)^m/(m+1)$ for
$z = 0.95$} \\ [1\jot] \hline \rule{0pt}{4\jot}%
$n$ & \multicolumn{1}{c}
{$\sum_{m=0}^{\infty} \frac{(-1)^{n+m} z^m}{n+m+2}$} &
\multicolumn{1}{c}{$z^{n+1} R_{\Ent {n/2}}^{(n - 2 \Ent {n/2})} (z)$ 
\rule[-10pt]{0pt}{25pt}} &
\multicolumn{1}{c}{$z^{n+1} 
r_{2 \Ent {n/2}}^{(n - 2 \Ent {n/2})} (z)$} &
\multicolumn{1}{c}{$z^{n+1} 
\mathcal{R}_{\Ent {n/3}}^{(n - 3 \Ent {n/3})} (z)$} \\ [1\jot]
 & & \multicolumn{1}{c}{Eq.\ (\ref{Ait_AccThrOrd})} &
\multicolumn{1}{c}{Eq.\ (\ref{Eps_AccThrOrd})} & 
\multicolumn{1}{c}{Eq.\ (\ref{tit_AccThrOrd})}
\\ [1\jot] \hline \rule{0pt}{4\jot}%
0  & $\phantom{-}0.312654\cdot 10^{0}$                 
& \multicolumn{1}{c}{$0$}                        
& \multicolumn{1}{c}{$0$} 
& \multicolumn{1}{c}{$0$} \\                      
1  & $-0.197206\cdot 10^{0}$		
& \multicolumn{1}{c}{$0$}                        
& \multicolumn{1}{c}{$0$}                       
& \multicolumn{1}{c}{$0$} \\                      
2  & $\phantom{-}0.143292\cdot 10^{0}$			
& $\phantom{-}0.620539\cdot 10^{-2}$    
& $\phantom{-}0.620539\cdot 10^{-2}$   
& \multicolumn{1}{c}{$0$} \\                      
3  & $-0.112324\cdot 10^{0}$		
& $-0.230919\cdot 10^{-2}$   
& $-0.230919\cdot 10^{-2}$  
& $\phantom{-}0.113587\cdot 10^{-2}$ \\  
4  & $\phantom{-}0.922904\cdot 10^{-1}$	
& $\phantom{-}0.109322\cdot 10^{-3}$    
& $\phantom{-}0.156975\cdot 10^{-3}$   
& $-0.367230\cdot 10^{-3}$ \\ 
5  & $-0.782908\cdot 10^{-1}$	
& $-0.333267\cdot 10^{-4}$   
& $-0.466090\cdot 10^{-4}$  
& $\phantom{-}0.148577\cdot 10^{-3}$ \\  
6  & $\phantom{-}0.679646\cdot 10^{-1}$	
& $\phantom{-}0.131240\cdot 10^{-5}$    
& $\phantom{-}0.413753\cdot 10^{-5}$   
& $\phantom{-}0.137543\cdot 10^{-5}$ \\  
7  & $-0.600373\cdot 10^{-1}$	
& $-0.371684\cdot 10^{-6}$   
& $-0.108095\cdot 10^{-5}$  
& $-0.392983\cdot 10^{-6}$ \\ 
8  & $\phantom{-}0.537619\cdot 10^{-1}$	
& $\phantom{-}0.111500\cdot 10^{-7}$    
& $\phantom{-}0.110743\cdot 10^{-6}$   
& $\phantom{-}0.131377\cdot 10^{-6}$ \\  
9  & $-0.486717\cdot 10^{-1}$	
& $-0.311899\cdot 10^{-8}$   
& $-0.266535\cdot 10^{-7}$  
& $\phantom{-}0.412451\cdot 10^{-9}$ \\  
10 & $\phantom{-}0.444604\cdot 10^{-1}$	
& $\phantom{-}0.689220\cdot 10^{-10}$   
& $\phantom{-}0.298638\cdot 10^{-8}$   
& $-0.139178\cdot 10^{-9}$ \\ 
11 & $-0.409189\cdot 10^{-1}$	
& $-0.199134\cdot 10^{-10}$  
& $-0.678908\cdot 10^{-9}$  
& $\phantom{-}0.475476\cdot 10^{-10}$ \\ 
12 & $\phantom{-}0.378992\cdot 10^{-1}$	
& $\phantom{-}0.282138\cdot 10^{-12}$   
& $\phantom{-}0.808737\cdot 10^{-10}$  
& $-0.316716\cdot 10^{-12}$ \\
[1\jot] \hline 
\end{tabular}
\end{table}

The zeros, which are found in columns 3 - 5 of Table I, occur because
Aitken's iterated $\Delta^2$ process and Wynn's epsilon algorithm can
only compute a rational approximant if at least three consecutive
partial sums are available, and because the iteration of Brezinski's
theta algorithm requires at least four partial sums.

The result in Table I show once more that Aitken's iterated $\Delta^2$
process and Brezinski's iterated theta algorithm are in the case of the
series (\ref{LogSer}) for $\ln(1+z)/z$ apparently more effective than
Wynn's epsilon algorithm.

The second group of results of this article -- (\ref{AitRemPhi}),
(\ref{EpsRem_phi}), and (\ref{titRemPsi}) and the corresponding
recursive schemes (\ref{RecPhi}), (\ref{rec_phi}), and (\ref{Rec_Psi})
-- can for example be used to demonstrate how rational approximants work
if a divergent power series is to be summed. 

Let us therefore assume that the partial sums, which occur in
(\ref{AitRemPhi}), (\ref{EpsRem_phi}), and (\ref{titRemPsi}), diverge if
the index becomes large. Then, a summation to a finite generalized limit
$f (z)$ can only be accomplished if $z^{n+2k+1} \Phi_{k}^{(n)} (z)$ and
$z^{n+2k+1} \varphi_{2k}^{(n)} (z)$ in (\ref{AitRemPhi}) and
(\ref{EpsRem_phi}), respectively, converge to the negative of $f_{n+2k}
(z)$, and if $z^{n+3k+1} \Psi_{k}^{(n)} (z)$ in (\ref{titRemPsi})
converges to the negative of $f_{n+3k} (z)$. 

Table II shows that this is indeed the case. We again consider the
infinite series (\ref{LogSer}) for $\ln(1+z)/z$, but this time we choose
$z = 5.0$, which is clearly outside the circle of convergence. We use
the numerical values of the partial sums $\sum_{m=0}^{n} (-z)^m/(m+1)$
with $0 \le n \le 10$ as input data in the recursive schemes
(\ref{RecPhi}), (\ref{rec_phi}), and (\ref{Rec_Psi}) to compute the
transformation terms in (\ref{AitRemPhi}), (\ref{EpsRem_phi}), and
(\ref{titRemPsi}). The transformation terms, which are listed in columns
3 - 5 of Table II, were chosen in agreement with (\ref{AitAppLim}),
(\ref{EpsAppLim}), and (\ref{TitAppLim}), respectively. All calculations 
were done using the floating point arithmetics of Maple.

\begin{table}[htb]
\begin{tabular}{lllll}
\multicolumn{5}{l}{\bf Table II: Convergence of Transformation Terms to 
the Partial Sums} \\
\multicolumn{5}{l}{Transformation of 
$\ln(1+z)/z = \sum_{m=0}^{\infty} (-z)^m/(m+1)$ for
$z = 5.0$} \\ [1\jot] \hline \rule{0pt}{4\jot}%
$n$ & \multicolumn{1}{c}{$\sum_{m=0}^{n} \frac{(-z)^m}{m+1}$} &
\multicolumn{1}{c}{$z^{n+1} \Phi_{\Ent {n/2}}^{(n - 2 \Ent {n/2})} (z)$ 
\rule[-10pt]{0pt}{25pt}} &
\multicolumn{1}{c}{$z^{n+1} 
\varphi_{2 \Ent {n/2}}^{(n - 2 \Ent {n/2})} (z)$} &
\multicolumn{1}{c}{$z^{n+1} 
\Psi_{\Ent {n/3}}^{(n - 3 \Ent {n/3})} (z)$} \\ [1\jot]
 & & \multicolumn{1}{c}{Eq.\ (\ref{AitRemPhi})} &
\multicolumn{1}{c}{Eq.\ (\ref{EpsRem_phi})} & 
\multicolumn{1}{c}{Eq.\ (\ref{titRemPsi})}
\\ [1\jot] \hline \rule{0pt}{4\jot}%
0  & $\phantom{-}0.1000000000\cdot 10^{1}$   
& \multicolumn{1}{c}{$0$} 
& \multicolumn{1}{c}{$0$} 
& \multicolumn{1}{c}{$0$}\\  
1  & $-0.1500000000\cdot 10^{1}$	 
& \multicolumn{1}{c}{$0$} 
& \multicolumn{1}{c}{$0$} 
& \multicolumn{1}{c}{$0$}\\  
2  & $\phantom{-}0.6833333333\cdot 10^{1}$	 
& $-0.6410256410\cdot 10^{1}$      
& $-0.6410256410\cdot 10^{1}$      
& \multicolumn{1}{c}{$0$}\\  
3  & $-0.2441666667\cdot 10^{2}$	 
& $\phantom{-}0.2467105263\cdot 10^{2}$       
& $\phantom{-}0.2467105263\cdot 10^{2}$       
& $\phantom{-}0.2480158730\cdot 10^{2}$\\        
4  & $\phantom{-}0.1005833333\cdot 10^{3}$	 
& $-0.1002174398\cdot 10^{3}$      
& $-0.1002155172\cdot 10^{3}$      
& $-0.1002604167\cdot 10^{3}$\\       
5  & $-0.4202500000\cdot 10^{3}$	 
& $\phantom{-}0.4205996885\cdot 10^{3}$       
& $\phantom{-}0.4205974843\cdot 10^{3}$       
& $\phantom{-}0.4206730769\cdot 10^{3}$\\        
6  & $\phantom{-}0.1811892857\cdot 10^{4}$	 
& $-0.1811533788\cdot 10^{4}$      
& $-0.1811532973\cdot 10^{4}$      
& $-0.1811533744\cdot 10^{4}$\\       
7  & $-0.7953732143\cdot 10^{4}$	 
& $\phantom{-}0.7954089807\cdot 10^{4}$       
& $\phantom{-}0.7954089068\cdot 10^{4}$       
& $\phantom{-}0.7954089765\cdot 10^{4}$\\        
8  & $\phantom{-}0.3544904563\cdot 10^{5}$	 
& $-0.3544868723\cdot 10^{5}$      
& $-0.3544868703\cdot 10^{5}$      
& $-0.3544868636\cdot 10^{5}$\\       
9  & $-0.1598634544\cdot 10^{6}$	 
& $\phantom{-}0.1598638127\cdot 10^{6}$       
& $\phantom{-}0.1598638125\cdot 10^{6}$       
& $\phantom{-}0.1598638127\cdot 10^{6}$\\        
10 & $\phantom{-}0.7279206365\cdot 10^{6}$	 
& $-0.7279202782\cdot 10^{6}$      
& $-0.7279202781\cdot 10^{6}$      
& $-0.7279202782\cdot 10^{6}$\\       
[1\jot] \hline 
\end{tabular}
\end{table}

The results in Table II show that a sequence transformation accomplishes
a summation of a divergent series by constructing approximations to the
actual remainders. Both the partial sums as well as the actual
remainders diverge individually if their indices become large, but the
linear combination of the partial sum and the remainder has a constant
and finite value for every index.

The fact, that the transformation terms in (\ref{AitRemPhi}),
(\ref{EpsRem_phi}), and (\ref{titRemPsi}) approach the negative of the
corresponding partial sums of course also implies that one should not
try to sum a divergent series in this way. The subtraction of two nearly
equal terms would inevitably lead to a serious loss of significant
digits.

In the next example, the transformation terms in (\ref{AitRemPhi}),
(\ref{EpsRem_phi}), and (\ref{titRemPsi}) will be used to make
predictions for unknown series coefficients. For that purpose, it is
recommendable to use a computer algebra system like Maple, and do all
calculations symbolically. If the coefficients of the series to be
transformed are exact rational numbers, the resulting rational
expressions are then computed exactly.

We use the symbolic expressions for the partial sums $\sum_{m=0}^{n}
(-z)^m/(m+1)$ with $0 \le n \le 12$ of the infinite series
(\ref{LogSer}) for $\ln(1+z)/z$ as input data in the recursive schemes
(\ref{RecPhi}), (\ref{rec_phi}), and (\ref{Rec_Psi}). The resulting
rational expressions $z^{13} \Phi_{6}^{(0)} (z)$, $z^{13}
\varphi_{12}^{(0)} (z)$, and $z^{13} \Psi_{4}^{(4)}$ with unspecified 
$z$ are then expanded, yielding predictions for the next series
coefficients that are exact rational numbers. Only in the final step,
the predictions for the next series coefficients are converted to
floating point numbers in order to improve readability:
\begin{subequations}
\label{f12_pred}
\begin{eqnarray}
\mathcal{A}_{6}^{(0)} & = & \sum_{m=0}^{12} \, \frac{(-z)^m}{m+1} 
\, - \, 0.07142857137 \, z^{13} \, + \, 0.06666666629 \, z^{14} 
\nonumber \\
& & \qquad
\, - \, 0.06249999856\,{z}^{15} \, + \, 0.05882352524\,{z}^{16}
\, + \, O \bigl( z^{17} \bigr) \, ,
\\
\epsilon_{12}^{(0)} & = & \sum_{m=0}^{12} \, \frac{(-z)^m}{m+1} 
\, - \, 0.07142854717 \, z^{13} \, + \, 0.06666649774 \, z^{14} 
\nonumber \\
& & \qquad
\, - \, 0.06249934843 \, {z}^{15} \, + \, 0.05882168762\,{z}^{16}
\, + \, O \bigl( z^{17} \bigr) \, ,
\\
\mathcal{J}_{4}^{(0)} & = & \sum_{m=0}^{12} \, \frac{(-z)^m}{m+1}
\, - \, 0.07142857148 \, {z}^{13} \, + \, 0.06666666684 \, {z}^{14}
\nonumber \\
& & \qquad
\, - \, 0.06249999986 \, {z}^{15} \, + \, 0.05882352708\,{z}^{16}
\, + \, O \bigl( z^{17} \bigr) \, ,
\\
\frac{\ln(1+z)}{z} & = & \sum_{m=0}^{12} \, \frac{(-z)^m}{m+1}
\, - \, 0.07142857143 \, {z}^{13} \, + \, 0.06666666667 \, {z}^{14}
\nonumber \\
& & \qquad
\, - \, 0.06250000000\,{z}^{15} \, + \, 0.05882352941\,{z}^{16}
\, + \, O \bigl( z^{17} \bigr) \, .
\end{eqnarray}
\end{subequations}
The accuracy of the prediction results in (\ref{f12_pred}) is quite
remarkable. The coefficients $\gamma_m = (-1)^m/(m+1)$ with $0 \le m \le
12$ are the only information that was used for the construction of the
transformation terms $z^{13} \Phi_{6}^{(0)} (z)$, $z^{13}
\varphi_{12}^{(0)} (z)$, and $z^{13} \Psi_{4}^{(0)}$, which were
expanded to yield the results in (\ref{f12_pred}). The accuracy of the
approximations to the next four coefficients should suffice for many
practical applications.

As in all other application, Wynn's epsilon algorithm is in
(\ref{f12_pred}) slightly but significantly less effective than
Aitken's iterated $\Delta^2$ process and Brezinski's iterated theta
algorithm.

Instead of computing the transformation terms $z^{13} \Phi_{6}^{(0)}
(z)$, $z^{13} \varphi_{12}^{(0)} (z)$, and $z^{13} \Psi_{4}^{(0)}$, it
is of course also possible to compute $\mathcal{A}_{6}^{(0)}$,
$\epsilon_{12}^{(0)}$, and $\mathcal{J}_{4}^{(0)}$ directly via their
defining recursive schemes, and to expand the resulting rational
expressions with a symbolic system like Maple. This would lead to the
same results. However, in order to extract the partial sum
$\sum_{m=0}^{12} (-z)^m/(m+1)$ from the rational approximants
$\mathcal{A}_{6}^{(0)}$, $\epsilon_{12}^{(0)}$, and
$\mathcal{J}_{4}^{(0)}$, one would have to compute their 12-th order
derivatives, and only the next derivatives would produce predictions to
unknown series coefficients. Thus, this approach can easily become very
expensive. In contrast, the use of the transformation terms requires
only low order derivatives of rational expressions.

If only the prediction of a single unknown term is to be done, then it
is of course much more efficient to use the recursive schemes
(\ref{RecG}), (\ref{Rec_g}), and (\ref{Rec_mclG}). The input data of
these recursive schemes are the coefficients of the series to be
transformed, and no differentiations have to be done.

\setcounter{equation}{0}

\section{Summary and Conclusions}
\label{Sec:SumCon}

As already mentioned in Section \ref{Sec:Intro}, it has become customary
in certain branches of theoretical physics to use Pad\'{e} approximants
to make predictions for the leading unknown coefficients of strongly
divergent perturbation expansions. This can be done by constructing
symbolic expressions for Pad\'{e} approximants from the known
coefficients of the perturbation series. A Taylor expansion of
sufficiently high order of such a Pad\'{e} approximants then produces
the predictions for the series coefficients which were not used for the
construction of the Pad\'{e} approximant. The Taylor expansion of the
symbolic expression can be done comparatively easily with the help of
powerful computer algebra systems like Maple or Mathematica, which are
now commercially available for a wide range of computers.

It is the purpose of this article to overcome two principal shortcomings
of the approach sketched above: Firstly, it is not necessary to rely
entirely on the symbolic capabilities of computers. Instead, it is
possible to construct recursive schemes, which either facilitate
considerably the symbolic tasks computers have to perform, or which
permit a straightforward computation of the prediction for the leading
unknown coefficient. Secondly, it is possible to use instead of Pad\'{e}
approximants other sequence transformations, as proposed by
Sidi and Levin \cite{SidLev83} and Brezinski \cite{Bre85}. It was shown
in \cite{We97} that this may lead to more accurate predictions. 

In this article, the prediction properties of Aitken's iterated
$\Delta^2$ process, Wynn's epsilon algorithm, and Brezinski's iterated
theta algorithm are studied.

As is well known \cite{Ba75,BaGM96}, a Pad\'{e} approximant can be
considered to be the solution of a system of linear equations for the
coefficients of its numerator and denominator polynomials. If this
system of linear equations has a solution, then it is automatically
guaranteed that the Pad\'{e} approximant satisfies the
accuracy-through-order relationship (\ref{PadeOrdEst}). In the case of
other sequence transformations, the situation is usually much more
difficult. They are usually not defined as solutions of systems of
linear equations, but via (complicated) nonlinear recursive schemes.

Since accuracy-through-order relationships of the type of
(\ref{PadeOrdEst}) play a very important role for the understanding of
the prediction properties of sequence transformations, it was necessary
to derive accuracy-through-order relationships for Aitken's iterated
$\Delta^2$ process, Wynn's epsilon algorithm, and Brezinski's iterated
theta algorithm on the basis of their defining recursive schemes.

Unfortunately, the defining recursive schemes (\ref{ItAit_1}),
(\ref{eps_al}), and (\ref{thetit_1}) are not suited for a construction
of accuracy-through-order relationships. They first had to be modified
appropriately, yielding the mathematically equivalent recursive schemes
(\ref{ItAit_2}), (\ref{CrossRule_2}), and (\ref{thetit_2}).

These alternative recursive schemes were the starting point for the
derivation of the accuracy-through-order relationships
(\ref{Ait_AccThrOrd}), (\ref{Eps_AccThrOrd}), and (\ref{tit_AccThrOrd})
and the corresponding recursive schemes (\ref{RecR}),
(\ref{Eps_AccThrOrd}), and (\ref{Rec_calR}) for the transformation error
terms. These relationships describe how the rational approximants
$\mathcal{A}_{k}^{(n)}$, $\epsilon_{2k}^{(n)}$, and
$\mathcal{J}_{k}^{(n)}$ differ from the function $f (z)$ which is to be
approximated.

With the help of these accuracy-through-order relationships, a second
group of results could be derived -- (\ref{AitRemPhi}),
(\ref{EpsRem_phi}), and (\ref{titRemPsi}) and the corresponding
recursive schemes (\ref{RecPhi}), (\ref{rec_phi}), and (\ref{Rec_Psi})
-- which describe how the rational approximants $\mathcal{A}_{k}^{(n)}$,
$\epsilon_{2k}^{(n)}$, and $\mathcal{J}_{k}^{(n)}$ differ from the
partial sums which were used for their construction. These differences
are expressed by the terms $z^{n+2k+1} \Phi_{k}^{(n)} (z)$, $z^{n+2k+1}
\varphi_{2k}^{(n)} (z)$, and $z^{n+3k+1} \Psi_{k}^{(n)} (z)$ which can
be computed via the recursive schemes (\ref{RecPhi}), (\ref{rec_phi}),
and (\ref{Rec_Psi}). 

The predictions for the leading unknown series coefficients can be
obtained by expanding symbolic expressions for these transformation
terms. The advantage of this approach is that the partial sums, which
are used for the construction of the rational approximants
$\mathcal{A}_{k}^{(n)}$, $\epsilon_{2k}^{(n)}$, and
$\mathcal{J}_{k}^{(n)}$ as well as of the transformation terms
$z^{n+2k+1} \Phi_{k}^{(n)} (z)$, $z^{n+2k+1} \varphi_{2k}^{(n)} (z)$,
and $z^{n+3k+1} \Psi_{k}^{(n)} (z)$, are already explicitly
separated. Consequently, only derivatives of low order have to be
computed. Moreover, the predictions for the leading unknown series
coefficient can be computed conveniently via the recursive schemes
(\ref{RecG}), (\ref{Rec_g}), and (\ref{Rec_mclG}). In this way, it is
neither necessary to construct symbolic expressions nor to differentiate
them.

Finally, in Section \ref{Sec:Applics} some applications of the new
results were presented.  In all applications of this article, Wynn's
epsilon algorithm was found to be less effective than Aitken's iterated
$\Delta^2$ process or Brezinski's iterated theta algorithm. Of course,
it remains to be seen whether this observation is specific for the
infinite series (\ref{LogSer}) for $\ln(1+z)/z$, which was used as the
test system, or whether it is actually more generally valid.
Nevertheless, the results presented in Section \ref{Sec:Applics} provide
further evidence that suitably chosen sequence transformations may
indeed be more effective than Pad\'{e} approximants. Consequently, one
should not assume that Pad\'{e} approximants produce by default the best
results in convergence acceleration and summation processes, and it may
well be worth while to investigate whether sequence transformations can
be found which are better adapted to the problem under consideration.

\section*{Acknowledgments}

My interest in Pad\'{e} approximants, sequence transformation,
convergence acceleration, and the summation of divergent series -- which
ultimately led to this article -- was aroused during a stay as a
Postdoctoral Fellow at the Faculty of Mathematics of the University of
Waterloo, Ontario, Canada. Special thanks to Prof.\ J.\
\v{C}\'{\i}\v{z}ek for his invitation to work with him, for numerous
later invitations to Waterloo, for his friendship, and the inspiring
atmosphere which he has been able to provide. Many thanks also to PD
Dr.\ H.\ Homeier for stimulating and fruitful discussions. Financial
support by the Fonds der Chemischen Industrie is gratefully
acknowledged.

\end{document}